\documentclass[10pt]{paper}%
\usepackage[T1]{fontenc}
\usepackage{amsmath,amssymb}
\usepackage{fullpage}
\usepackage{amsthm}
\usepackage{multirow}
\usepackage{xcolor}
\usepackage{pifont}
\usepackage{graphicx}
\usepackage[all]{xy}
\usepackage{amsfonts}
\usepackage{thmtools}
\usepackage{hyperref}
\usepackage[nameinlink]{cleveref}
\usepackage{enumitem}

\providecommand{\U}[1]{\protect\rule{.1in}{.1in}}
\providecommand{\U}[1]{\protect\rule{.1in}{.1in}}
\providecommand{\U}[1]{\protect\rule{.1in}{.1in}}
\providecommand{\U}[1]{\protect\rule{.1in}{.1in}}
\providecommand{\U}[1]{\protect\rule{.1in}{.1in}}
\newcommand{\ulambda}{{\boldsymbol{\lambda}}}

\DeclareMathOperator{\id}{Id}
\DeclareMathOperator{\seq}{seq}
\DeclareMathOperator{\st}{st}

\newtheorem{Th}{Theorem}[section] 
\newtheorem{lemma}[Th]{Lemma}

\newtheorem{Prop}[Th]{Proposition}

\theoremstyle{remark}
\newtheorem{Rem}[Th]{Remark}{\rmfamily}
\theoremstyle{definition}
\newtheorem{Def}[Th]{Definition}{\rmfamily}
\newtheorem{exa}[Th]{Example}{\rmfamily}

\newcommand\blfootnote[1]{%
  \begingroup
  \renewcommand\thefootnote{}\footnote{#1}%
  \addtocounter{footnote}{-1}%
  \endgroup
}

\def\arxivprefix{https://arxiv.org}
\newcommand*\arxiv[1]{\href{\arxivprefix/abs/#1}{arXiv:#1}}

\begin{document}

\title{On the computation of  the canonical basis for irreducible highest weight $U_q(\mathfrak{gl}_{\infty})$-modules}
\author{Nicolas Jacon and Abel Lacabanne} 
\maketitle
\date{}
\blfootnote{\textup{2020} \textit{Mathematics Subject Classification}: \textup{05E10,17B37,20C08}} 
\begin{abstract}
We study canonical basis elements in higher-level Fock spaces associated with the quantum group $U_q(\mathfrak{gl}_\infty)$, which are conjecturally related to Calogero--Moser theory for complex reflection groups. We generalize the Leclerc--Miyachi formula to arbitrary levels by introducing new explicit constructions based on symbols, including a column removal theorem and closed formulas in several cases. These results provide explicit descriptions of canonical basis elements with applications to Calogero--Moser cellular characters and to the decomposition matrices of Ariki--Koike algebras.
\end{abstract}


\section{Introduction}
An important problem in Lie theory is to generalize Kazhdan--Lusztig theory to the class of complex reflection groups. This theory originates with the construction of a distinguished basis known as the Kazhdan--Lusztig basis for the Hecke algebra of a Coxeter group, and it yields a partition of the elements of the Weyl group into left, right, and two-sided cells.

Such a basis does not exist for the entire class of complex reflection groups, and in particular for the class of imprimitive complex reflection groups, namely wreath products of a cyclic group with a symmetric group. However, more recently, Bonnaf\'e and Rouquier \cite{BR} proposed a generalization of the notion of cells via the so-called Calogero--Moser cells, defined through the study of the Calogero--Moser space associated with a complex reflection group. These cells are conjectured to play a role analogous to that of Kazhdan--Lusztig cells for Coxeter groups. This correspondence has been established in several cases, notably for the symmetric group \cite{BGW} and for dihedral groups with equal parameters \cite{BG}.

In the same spirit, a natural generalization of constructible characters for Weyl groups is given by the so-called Calogero--Moser cellular characters. In \cite{La}, the second author conjectured that these Calogero--Moser characters in type $G(l,1,n)$ can be computed via the expansion of the canonical basis in an irreducible highest weight $U_q(\mathfrak{gl}_{\infty})$-module of level~$l$. These canonical basis elements depend on the choice of an $l$-tuple of integers $\mathbf{v}=(v_1,\dots,v_l)$, called a multicharge, which parametrizes the highest weight. Each canonical basis element is associated with a certain combinatorial object, known as a (standard) symbol, arising from Kashiwara's crystal basis theory.

{It turns out that these canonical basis elements can be computed explicitly in certain cases. This is notably true for $l=2$, where Leclerc and Miyachi \cite{HL} provided an explicit and elegant formula for the canonical basis elements. Moreover, at this level, they showed that the constructible characters in type $B_n$ correspond to canonical basis elements of level~$2$. More generally, while algorithms exist to compute these canonical basis elements \cite{LT,Jal}, they do not provide explicit formulas.}

The aim of this paper is to obtain new closed and explicit formulas for canonical basis elements at higher levels. By Ariki's theorem, these elements also describe the decomposition matrices of Ariki--Koike algebras at a {nonroot} of unity. In \cite{JL}, we presented a first systematic study of these canonical basis elements and their relationship with Calogero--Moser theory. In the present work, we pursue this study further by generalizing the Leclerc--Miyachi formula and by providing the first closed formulas for canonical basis elements at level $l>2$. Our main results are as follows:

\begin{itemize}
\item In \Cref{asy}, we study certain special multicharges that may be viewed as a generalization of the asymptotic case previously investigated in \cite{La}.

\item In \Cref{removal}, we establish an analogue of the column removal theorem for canonical basis elements using the theory of symbols. We show that, in certain cases, the canonical basis element associated with a symbol $S$ can be obtained from that associated with another symbol derived from $S$ by removing specific entries.

\item In \Cref{ordered}, we provide an explicit formula for certain canonical basis elements at arbitrary levels and for arbitrary multicharges. In particular, for each block, there exists a distinguished symbol that can be computed explicitly using our formula.

\item Finally, in \Cref{thm:main_formula}, we give a formula for canonical basis elements at higher levels, which naturally generalizes the Leclerc--Miyachi formula for the level~$2$ case.
\end{itemize}

The paper is organized as follows. In \Cref{sec:algebra-sl}, we recall the representation theory and combinatorics of the Fock space, including the notion of a symbol associated with a multipartition, which plays a central role throughout the paper. In \Cref{removal}, we establish the column removal theorem and extend the asymptotic case previously studied by the second author. Finally, the last section is divided into two subsections, each devoted to presenting and proving the main results concerning the structure of canonical basis elements.

\medskip

{\bf Acknowledgements.} The authors would like to thank the referee for pointing a mistake in a previous version of this paper. A.L. also would like to thank Dani Tubbenhauer for useful discussions on the level $3$ case.

\section{The algebra $U_q(\mathfrak{gl}_{\infty})$ and the Fock space $\mathcal{F}_{\mathbf{v}}$}
\label{sec:algebra-sl}

We here define our main objects of study: the quantum group $U_q(\mathfrak{gl}_{\infty})$ and its action on the Fock space. After several combinatorial definitions, we introduce the notion  of  canonical basis. Our main references here are \cite{GJ} and \cite{Alivre}.

\subsection{Definitions}
We consider the \emph{quantum group} $U_q(\mathfrak{gl}_{\infty})$. This is the  $\mathbb{Q}(q)$-algebra with generators $E_i,F_i$ and $L_i^{\pm 1}$ for $i \in \mathbb{Z}$ subject to the following relations:
\[
  L_iL_i^{-1} = 1 = L_i^{-1}L_i,\quad L_i L_j = L_j L_i, \quad L_iE_j = q^{\delta_{i,j}-\delta_{i,j+1}}E_jL_i,
\]
\[
  L_iF_j = q^{-\delta_{i,j}+\delta_{i,j+1}}F_jL_i, \quad [E_i,F_j] = \delta_{i,j}\frac{L_iL_{i+1}^{-1}-L_i^{-1}L_{i+1}}{q-q^{-1}}
\]
together with  the Serre relations
\[
  E_i^2E_j - [2] E_iE_jE_i + E_jE_i^2 = 0,\  F_i^2F_j - [2] F_iF_jF_i + F_jF_i^2 = 0,\text{ if }|i-j| = 1,
\]
\[
  [E_i,E_j] = 0,\  [F_i,F_j]=0, \text{ if }|i-j| > 1,
\]
where for $r \in \mathbb{Z}_{\geq 0}$, we let $[r]=\displaystyle{\frac{q^r-q^{-r}}{q-q^{-1}}}$ and $[r]!=[r][r-1]\ldots [1]$. 
We define the divided powers $E_i^{(r)}$ and $F_i^{(r)}$ by
\[
  E_i^{(r)}=\frac{E_i}{[r]!}\quad\text{and}\quad F_i^{(r)}=\frac{F_i}{[r]!}.
\]
The \emph{fundamental weights} of $U_q(\mathfrak{gl}_{\infty})$ are denoted by $(\Lambda_k)_{k\in\mathbb{Z}}$. We now fix {$l > 0$} and $\mathbf{v}=(v_1,\ldots,v_l)\in\mathbb{Z}^l$. We want to show how one can realize the irreducible highest weight module with weight $\Lambda_{v_1}+\ldots +\Lambda_{v_l}$ (we say that such a module is of level $l$) as submodule of the Fock space, defined in the subsections below. We refer to \cite{Alivre,GJ,Jlivre} for more details on the definition below.

We also endow $U_q(\mathfrak{gl}_{\infty})$ with a structure of a Hopf algebra, whose coproduct $\Delta$, antipode $S$ and counit $\varepsilon$ are given by
\begin{align*}
  \Delta(E_i) &= E_i\otimes 1 + L_i^{-1}L_{i+1}\otimes E_i, & S(E_i) &= -L_iL_{i+1}^{-1}E_i, & \varepsilon(E_i) &= 0,\\
  \Delta(F_i) &= F_i\otimes L_iL_{i+1}^{-1} + 1\otimes F_i,& S(F_i) &= -F_iL_{i}^{-1}L_{i+1}, & \varepsilon(F_i) &= 0,\\
  \Delta(L_i) &= L_i\otimes L_i ,& S(L_i) &=L_i^{-1}, & \varepsilon(L_i) &= 1.
\end{align*}

\subsection{Multipartitions and symbols}

A \emph{partition} $\lambda = (\lambda_1,\lambda_2,\ldots)$ of an integer $n$ is a {nonincreasing} sequence $\lambda_1 \geq \lambda_2 \geq \cdots $ of {nonnegative} integers of sum  $n$. Such a sequence is eventually zero and we may only write the {nonzero} terms. If $\lambda=(\lambda_1,\lambda_2,\ldots)$ is a partition of an integer, we denote by $\lvert \lambda \rvert$ the sum $\sum_{i\geq 1}\lambda_i$, the \emph{size} of the partition. 

Let $l \in \mathbb{N}$. An $l$-partition (or multipartition)  $\ulambda = (\lambda^{(1)},\ldots,\lambda^{(l)})$ of $n$ is an $l$-tuple of partitions such that $\sum_{i=1}^l \lvert\lambda^{(i)}\rvert = n$, which is called the size of $\ulambda$. 
We fix  $\mathbf{v}=(v_1,\ldots,v_l)$  such that 
\[
  \mathbf{v} \in \mathcal{A}^l=\{ \mathbf{v}=(v_1,\ldots,v_l)\in \mathbb{Z}^l\ |\ v_1\geq v_2 \geq \ldots \geq v_l \}.
\]
This is called a \emph{multicharge}. 

The \emph{Fock space} has a basis indexed by all $l$-partitions, but the action of the quantum enveloping algebra $U_q(\mathfrak{gl}_{\infty})$ is better understood in terms of $l$-symbols that we now define. For $s\in \mathbb{Z}$, we denote by $\mathfrak{B}_{s}$ the subset of $\mathbb{Z}^{\mathbb{Z}_{\leq s}}$ given by  the sequences $(\beta_j)_{j \leq s}$ such that $\beta_{j-1} < \beta_j$ and $\beta_j = j$ for $j \ll s$. This is the set of \emph{$\beta$-numbers}. 

An \emph{$l$-symbol} (or symbol if $l$ is understood) of \emph{multicharge} $\mathbf{v}$ is an $l$-tuple $(\beta^1,\ldots,\beta^l)$ where $\beta^i=(\beta^i_j)_{j \leq v_i}\in \mathfrak{B}_{v_i} $. We represent such a symbol as:
\[
  S =
  \begin{pmatrix}
    \beta^1 \\
    \beta^{2} \\
    \vdots \\
    \beta^l
  \end{pmatrix}
  =
  \begin{pmatrix}
    \cdots & \beta^{1}_{v_l-1} & \beta^{1}_{v_l} & \beta^{1}_{v_l+1} & \cdots & \beta^{1}_{v_{2}-1} &  \beta^{1}_{v_{2}} & \beta^{1}_{v_{2}+1} & \cdots & \beta^{1}_{v_1} \\
    \cdots & \beta^{2}_{v_l-1} & \beta^{2}_{v_l} & \beta^{2}_{v_l+1} & \cdots & \beta^{2}_{v_{2}-1} &  \beta^{2}_{v_{2}} &  &  &  \\
    \vdots & \vdots & \vdots & \vdots & \vdots & \vdots & \vdots & & & \\
    \cdots & \beta^{l}_{v_l-1} & \beta^{l}_{v_{l}} & & & & & & & 
  \end{pmatrix}.
\]

The size of such a symbol is the integer $\sum_{i=1}^l \sum_{j \leq v_i}(\beta^i_{j}-j)$. Since $\beta^i_j = j$ for $j \ll 0$, the size of a symbol is well defined. If this size is $0$, we say that the symbol is empty and we denote it by $\emptyset_\mathbf{v}$. This is the case when $\beta^{i}_{j} = j$ for all relevant $i$ and $j$. The set of symbols associated with the multicharge $\mathbf{v}$ is denoted by $\mathcal{S}_\mathbf{v}$ and the set of elements in $\mathcal{S}_\mathbf{v}$ of size $n$ by $\mathcal{S}_\mathbf{v}(n)$. 

Multipartitions and symbols of a given multicharge are in bijection. Given an $l$-partition $\ulambda =(\lambda^{(1)},\ldots,\lambda^{(l)})$ of $n$, we associate the $l$-symbol $S_{\lambda}=(\beta^1,\ldots,\beta^l)$ of multicharge $\mathbf{v}$ and size $n$ defined by {$\beta^i_j = \lambda^{(i)}_{v_i-j+1}+j$} for $1 \leq i \leq l$ and $j \leq v_i$. Conversely, given an $l$-symbol $S$ of multicharge $\mathbf{v}$ and size $n$, we associate the $l$-partition $\ulambda_S =(\lambda^{(1)},\ldots,\lambda^{(l)})$ of $n$ such that $\lambda^{(i)}_j = \beta^{i}_{v_i-j+1}-v_i+j-1$. These two operations are inverse to each other.

We say that a symbol $S=(\beta_1,\ldots,\beta_l)$ is {\it standard} if $\beta^{i}_{j} \leq \beta^{i+1}_{j}$ for all $1 \leq i < l$ and for all  $j\leq v_{i}$. The set of standard symbols is denoted by $\mathcal{S}^{\st}_\mathbf{v}$ and the set of standard symbols of size $n$ is denoted by $\mathcal{S}^{\st}_\mathbf{v}(n)$.

\begin{exa}
The $4$-symbol 
\[
  S=\left( 
    \begin{array}{ccccccc}
      \ldots & 0 & 1 & 2 & 4 & 5 & 6 \\
      \ldots &0 & 1 & 3 & 5  & 7 & 8\\
      \ldots &0 & 1 & 4 \\
      \ldots &0 & 2& 4 \\
    \end{array}\right)
\]
is a standard symbol. The associated multicharge is $(5,5,2,2)$ and the associated multipartition is the $4$-partition  $((1,1,1),(3,3,2,1),(2),(2,1))$. The size is $17$.  
\end{exa}

\begin{exa}
Consider the $3$-symbol:
\[
  S=\left( 
    \begin{array}{cccccc}
      \ldots & 0 & 1 & 3 & 5 & 7 \\
      \ldots &0 & 1 & 3 & 5  \\
      \ldots &0 & 1 & 4  & 6\\
    \end{array}\right)
\]
The associated multicharge is $(4,3,3)$ and the associated multipartition is $((3,2,1),(2,1),(3,2))$. The size is $14$.
\end{exa}

\subsection{Action on the Fock space}

Given $s\in \mathbb{Z}$,  we now construct the $U_q(\mathfrak{gl}_{\infty})$-module $V(\Lambda_s)$ of  highest weight the fundamental weight $\Lambda_s$. This  has a $\mathbb{Q}(q)$-basis  given  by sequences $\beta=(\beta_j)_{j\leq s}$ in  $\mathfrak{B}_{s}$. The action of the generators is given as follows:
\[
  E_i\cdot \beta =
  \begin{cases}
    (\beta\setminus \{i+1\}) \cup \{i\} & \text{if } i+1 \in \beta \text{ and } i\not\in \beta,\\
    0 & \text{otherwise},
  \end{cases}
  \quad
  F_i\cdot \beta =
  \begin{cases}
    (\beta\setminus \{i\}) \cup \{i+1\} & \text{if } i \in \beta \text{ and } i+1\not\in \beta,\\
     0 & \text{otherwise},
  \end{cases}
\]
\[    
  L_i\cdot  {\beta} =
  \begin{cases}
    q \beta & \text{if } i \in \beta \text{ and } i+1\not\in \beta,\\
    q^{-1} \beta & \text{if } i+1 \in \beta \text{ and } i\not\in \beta,\\
    \beta & \text{otherwise}.
  \end{cases}
\]

\begin{Def}
  Let $\mathbf{v}=(v_1,\ldots,v_l)$ be a multicharge with $v_1 \geq \cdots \geq v_l$. The \emph{Fock space} $\mathcal{F}_{\mathbf{v}}$ is the tensor product $V(\Lambda_{v_1})\otimes \cdots \otimes V(\Lambda_{v_l})$, where the action of $U_q(\mathfrak{gl}_{\infty})$ is given by the coproduct $\Delta$.
\end{Def}

Therefore, the Fock space $\mathcal{F}_{\mathbf{v}}$ has a basis indexed by the $l$-symbols of multicharge $\mathbf{v}$, or equivalently by the $l$-partitions of integers.

The $U_q(\mathfrak{gl}_{\infty})$-subspace $V_{\mathbf{v}}$ generated by the vector $\emptyset_{\mathbf{v}}$ is isomorphic to the representation $V\left(\Lambda_{v_1}+\cdots+\Lambda_{v_l}\right)$ of highest weight $\Lambda_{v_1}+\cdots+\Lambda_{v_l}$. From the definition of the divided powers, we get the following classical result (see \cite[Prop. 6.2.7]{GJ}).

\begin{Prop}\label{divided}
Let $S=(\beta^1,\ldots,\beta^l)$ be a symbol. Let $i\in \mathbb{Z}$ and let $a\in \mathbb{N}$. Then we have:
\[
  E_i^{(a)}\cdot {S} =\sum_{S'} q^{-N^1(S,S')} {S'}
\]
where the sum is taken over all the symbol  $S'=(\gamma^1,\ldots,\gamma^l)$ obtained from $S$ by replacing exactly $a$ entries $i+1$ with $i$, say in rows $j_1, \ldots, j_a$, and where
\[
  N^1 (S,S')=\sum_{1\leq r \leq a} 
  \left( \sharp\{k<j_r \mid i\in \gamma^k\}-\sharp\{k<j_r\mid i+1\in \beta^k\}\right).
\]
Similarly, we have:
\[
  F_i^{(a)}\cdot {S} =\sum_{S'} q^{N^2(S,S')} {S'},
\]
where the sum is taken over all the symbol $S'=(\gamma^1,\ldots,\gamma^l)$ obtained from $S$ by replacing exactly $a$ entries $i$ with $i+1$, say in rows $j_1, \ldots, j_a$, and where
\[
  N^2(S,S')=\sum_{1\leq r \leq a} 
  \left(\sharp\{k>j_r \mid i\in \gamma^k\}-\sharp\{k>j_r \mid i+1\in \beta^k\}\right).
\]
\end{Prop}

\subsection{Monomial and canonical bases}

We now define the Kashiwara operators $\widetilde{F}_i$ and $\widetilde{E}_i$. To do this, we define the $i$-signature of a symbol $S$ as follows. We read the symbol from top to bottom and construct a word as follows: we write  $+$ if we meet a value $i+1$ and $-$ if we meet the value $i$ (if $i$ and $i+1$ are in the same row, we write an occurence $+-$) . We then delete repeatedly all the occurences of $+-$ in the word. We obtain a sequence of $N_1$ ``$-$'' and then $N_2$  ``$+$''. 
 \begin{enumerate}
\item If $N_1=0$ we let  $\widetilde{F}_iS=0$. Otherwise, the rightmost ``$-$'' in the  word correspond to an element $i$. The symbol  $\widetilde{F}_iS$ is then obtained  from $S$ by replacing this $i$ with $i+1$. 
\item  If  $N_2=0$ we let  $\widetilde{E}_iS=0$. Otherwise, the leftmost ``+'' in the  word correspond to an element $i+1$. the symbol  $\widetilde{E}_iS$ is then obtained  from $S$ by replacing this $i+1$ with $i$. 
\item We also define $\varepsilon_i (S)=N_1$ and $\varphi_i (S)=N_2$. 
\end{enumerate}

Note that we have:
\[
  \varepsilon_i(T)=\max \left\{ N\geq 0 \ \middle\vert\ \widetilde{E}_i^N T\neq 0\right\}
\]
and 
\[
  \varphi_i(T)=\max \left\{N\geq 0 \ \middle\vert\ \widetilde{F}_i^N  T\neq 0\right\}.
\]
By Kashiwara and Lusztig theory (see \cite{Kashi}), we get:

\begin{Prop}\label{Kas}
  Assume that $\mathbf{v}\in \mathcal{A}^l$. Let $S$ be a symbol of size $n$. Then  there exists $(i_1,\ldots,i_n)$ such that 
  \[
    \widetilde{F}_{i_1}\ldots \widetilde{F}_{i_n}\emptyset_\mathbf{v}  =S
  \]
  if and only if $S\in \mathcal{S}^{\st}_{\mathbf{v}}$.
\end{Prop}

\begin{exa}
  Assume that $\mathbf{v}=(3,2,1) $. We consider the following $3$-symbol:
  \[
    S=
    \left(
      \begin{array}{ccccc}
        \ldots & 0 & 1 & 3 & 4 \\
        \ldots & 0 & 1 & 4 & \\
        \ldots & 0 & 2 & &
      \end{array}
    \right)
  \]
  This is a standard symbol. We have $\widetilde{F}_{1}\widetilde{F}_{2}\widetilde{F}_{3}^2\widetilde{F}_{2}\emptyset_\mathbf{v} = S$.
\end{exa}

Let $x\mapsto \overline{x}$ be the $\mathbb{Q}$-linear algebra involution of $U_q(\mathfrak{gl}_{\infty})$ given by
\begin{align*}
  \overline{q}&=q^{-1}, & \overline{L_i}&= L_i^{-1}, & \overline{E_i}&=E_i, & \overline{F_i} = F_i.
\end{align*}
Any element $v\in V_{\mathbf{v}} $ can then be expressed as $v=x \cdot  \emptyset_\mathbf{v} $ for some $x\in U_q(\mathfrak{gl}_{\infty})$ and we set $\overline{v} = \overline{x} \cdot  \emptyset_\mathbf{v} $. This gives a well-defined involution on $V_{\mathbf{v}}$. 
Let $R$ be the subring of $\mathbb{Q}(q)$ of rational functions regular at $q=0$ and consider $\mathcal{F}_{\mathbf{v},R}$ the $R$-sublattice of $\mathcal{F}_{\mathbf{v}}$ generated by the symbols. By \cite{L90} we have:

\begin{Th}
  There exists a unique basis $\left\{G(S)\ \middle\vert\ \ S\in \mathcal{S}^{\st}_\mathbf{v} \right\}$ of $V_{\mathbf{v}}$ such that
  \begin{itemize}
  \item $\overline{G(S)} = G(S)$,
  \item $G(S) = S  \mod q\mathcal{F}_{\mathbf{v},R}$.
  \end{itemize}
\end{Th}

The basis $(G(S))$ of $V_{\mathbf{v}}$ is the \emph{canonical basis} of $V_{\mathbf{v}}$. It is easy to see that if 
\[
  G(S)=\sum_{T\in \mathcal{S}_\mathbf{v}} \alpha_{S,T} (q) T,
\]
then if $\alpha_{S,T} (q)\neq 0$, then the multiset of elements appearing in $T$ is the same as the one appearing in $S$.
In fact, one can define an equivalence relations on the set of symbols as follows: we say that $T$ and $T'$ are equivalent if and only if there exist a standard symbol $S$ such that  $\alpha_{S,T} (q) \neq 0$ and $\alpha_{S,T'} (q) \neq 0$. 
Our equivalence relation is the transitive closure of this relation. We thus have a partition  of the set of symbols in these equivalence classes, which are called the blocks. It is known that two symbols  $T$ and $T'$ are in the same block if and only if the multiset of elements appearing in $T$ is the same as the one appearing in $T'$, which implies that the blocks are also the Calogero--Moser families .

\begin{Rem}\label{plus}
Assume that $\mathbf{v}=(v_1,\ldots,v_l)$ is such that  $v_1 \geq \cdots \geq v_l$ and denote $\mathbf{v}^{+}=(v_1+1,\ldots,v_l+1)$. If $T$ is a symbol then denote $T^{+}$ the symbol obtained by adding $1$ to all the entries of the symbol. If  $S\in \mathcal{S}^{\st}_{\mathbf{v}}$ and if we have
\[
  G(S)=\sum_{T\in \mathcal{S}_\mathbf{v}} \alpha_{S,T}(q) T,
\]
then it is trivial to see that $S^+\in \mathcal{S}^{\st}_{\mathbf{v}^+} $ and 
\[
  G(S^+)=\sum_{T\in \mathcal{S}_\mathbf{v}} \alpha_{S,T}(q) T^+.
\]
\end{Rem}

Let $(i_1,\ldots,i_r)\in \mathbb{Z}^r$  be such that $i_j\neq i_{j+1}$ for $j=1,\ldots, r-1$ and let $(a_1,\ldots, a_r)\in \mathbb{N}^r$. Then we denote 
\[
  \underline{i}=(\underbrace{i_1,\ldots,i_1}_{a_1},\ldots,\underbrace{i_r,\ldots,i_r}_{a_r})
\]
and
\[
  F_{\underline{i}} \cdot \emptyset_{\mathbf{v}} = F^{(a_1)}_{i_1}\ldots F^{(a_r)}_{i_r} \cdot \emptyset_{\mathbf{v}}.
\]
We say that $G(S)$ is monomial if there exists $\underline{i}$ as above such that 
\[
  G(S)=F_{\underline{i}} \cdot \emptyset_{\mathbf{v}}.
\]

In the case $l=2$, Leclerc and Miyachi \cite{HL} have given a closed formula for the canonical basis  elements that we now recall. Let $(\gamma,\beta)\in \mathfrak{B}_t \times \mathfrak{B}_s$ with $t\geq s$ {such that the symbol $(\gamma,\beta)$ is standard}. Then, we define an injection 
\begin{align}\label{Lmy}
\Psi \colon  \beta \to \gamma
\end{align}
{as follows.} There exists {$j\in \mathbb{Z}$} such that $\beta_k=\gamma_k=k$ for all $k\leq j$. We then define $\Psi (\beta_k)=\gamma_k$ for all $k\leq j$. For $k>j$, we inductively  define:
\[
  \Psi (\beta_k)=\max \left\{\gamma_j \ \middle\vert\ \gamma_j \leq \beta_k,\text{ and }\forall r<k,\  \gamma_j\neq \Psi (\beta_r)\right\}.
\]
Let $S=(\beta^1,\beta^2)$ be a {standard} $2$-symbol and consider the injection $\Psi \colon \beta^{2} \to \beta^{1}$ as above. The pairs $(\Psi(\beta^2_j),\beta^2_j)$ with $\Psi(\beta^2_j)\neq \beta^2_j$ are called the {\it pairs} of $S$. We denote by $C(S)$ the set of symbols obtained from $S$ by permuting some pairs in $S$ and reordering the rows. For $T\in C(S)$, the number $n(T)$ is the number of pairs which have been permuted from $S$ to obtain $T$.

\begin{Th}[\cite{HL}]
  Assume that $l=2$. For all standard symbol $S$, we have 
  \[
    G(S)= \sum_{T \in C(S)} q^{n(T)}T.
  \]
\end{Th}

\begin{exa}
  Take $\mathbf{v}=(1,0)$ and the following standard symbol:
  \[
    S=
    \left(
      \begin{array}{ccccc}
        \ldots & 0 & 1 & 3 & 5 \\
        \ldots & 0 & 2 & 7 &  
      \end{array}
    \right)
  \]
  Then the pairs are $(1,2)$ and $(5,7)$. As a consequence, we get:
  \[
    G(S)=
    \left(
      \begin{array}{ccccc}
        \ldots & 0 & 1 & 3 & 5 \\
        \ldots & 0 & 2 & 7 & 
      \end{array}
    \right)
    +
    q\left(
      \begin{array}{ccccc}
        \ldots & 0 & 2 & 3 & 5 \\
        \ldots & 0 & 1 & 7 & 
      \end{array}
    \right)
    +
    q\left(
      \begin{array}{ccccc}
        \ldots & 0 & 1 & 3 & 7 \\
        \ldots & 0 & 2 & 5 & 
      \end{array}
    \right)
    +q^2
    \left(
      \begin{array}{ccccc}
        \ldots & 0 & 2 & 3 & 7 \\
        \ldots & 0 & {1} & 5 &  
      \end{array}
    \right).
  \]
\end{exa} 

\section{Removal Theorem}

From the datum of certain canonical basis elements, we here present two results which allow the computation of new canonical basis ones. 

\subsection{Asymptotic cases}\label{asy}

In this subsection, we fix $n\in \mathbb{N}$. We still assume that we have $v_1 \geq v_2 \geq \ldots \geq v_l$. In addition, we assume that there exist integers {$k_2<k_3 <\ldots < k_{r}$} in $\{2,\ldots,l\}$ such that 
 \[
   v_{k_{j}-1}-v_{k_j} \geq n 
 \]
 for all $j=1,\ldots,r-1$. If $r=l$ and if {$(k_2,\ldots,k_{r})=(2,\ldots,l)$}, this is the ``asymptotic case'' studied in \cite[\S 2]{La}. We denote for $j=1,\ldots,r$
 \[
   \mathbf{v}^j=(v_{k_{j}},\ldots,v_{k_{j+1}-1}),
 \]
 with the convention that {$k_1=1$} and {$k_{r+1}=l+1$}. We also let $l_j=k_{j+1}-k_{j}$ so that $\sum_{1\leq i \leq r} l_i=l$.

Fix $j\in \{1,\ldots,r\}$. We have the set of canonical basis elements for $V_{\mathbf{v}^j}$
\[
  \left\{ G(S)\ \middle\vert\ S \in \mathcal{S}^{\st}_{\mathbf{v}^j}(m)\right\}
\]
in the Fock space with $m\leq n$. For each $S \in \mathcal{S}^{\st}_{\mathbf{v}^j}(m)$, we have 
\[
  G(S) = \sum_{\underline{i} \in J(S)}  a^S_{\underline{i}}(q) F_{\underline{i}} \cdot \emptyset_{\mathbf{v}^j} = \sum_{T\in \mathcal{S}_{\mathbf{v}^j}(m)} \alpha^j_{S,T}(q) T,
\]
with $\alpha^j_{S,T}(q)\in \mathbb{Z}[q]$ and $a^S_{\underline{i}}(q) \in \mathbb{Z}[q,q^{-1}]$. Here  $J(S)$ is a certain subset of $\mathbb{Z}^{m}$ and 
 $a^S_{\underline{i}}(q)$  is bar invariant. In the following, we will 
  mostly write  $a_{\underline{i}}(q)$ instead of  $a^S_{\underline{i}}(q)$ to keep the notation light (if there is no risk of confusion).
  The following proposition then follows directly from the definition of standard symbols (see \cite[5.5.18]{GJ}):
  
\begin{Prop}\label{prop:standard-asym}
  Let $S$ be a symbol of size $n$ associated with the multicharge $\mathbf{v}$. Then $S$ is standard if and only if there exists $n_1,\ldots,n_r\in \mathbb{N}$ with $n=n_1+\dots+n_r$ and $S_j \in \mathcal{S}^{\st}_{\mathbf{v}^j}(n_j)$ such that $S = (S_1,\ldots,S_r)$.

  In other words, we can split a standard symbol $S$ in $r$ smaller standard symbols $S_1,\dots,S_r$.
\end{Prop}

We then obtain a formula for the canonical basis $G(S)$ with $S=(S_1,\ldots,S_r)$ a standard symbol of size $n$ from the canonical basis elements $G(S_1),\dots,G(S_r)$.

\begin{Prop}
  Let $S$ be a standard symbol of size $n$. Write $S=(S_1,\dots,S_r)$ as in \Cref{prop:standard-asym}, with $S_j \in \mathcal{S}^{\st}_{\mathbf{v}^j}(n_j)$ and $n_1+\ldots+n_r=n$. Then
  \begin{align*}
    G(S)&=\sum_{(\underline{i}_1,\ldots,\underline{i}_r)\in J(S_1) \times \ldots \times J(S_r)}
    a_{\underline{i}_1}(q)\ldots a_{\underline{i}_r}(q) F_{\underline{i}_1} \ldots F_{\underline{i}_r}\cdot \emptyset_{\mathbf{v}}\\
    &= \sum_{(T_1,\ldots,T_r)\in \mathcal{S}_{\mathbf{v}^1}(n_1) \times \dots \times \mathcal{S}_{\mathbf{v}^r}(n_r)} \alpha^1_{S_1,T_1}(q)\ldots \alpha^r_{S_r,T_r}(q) (T_1,\ldots,T_r).
  \end{align*}
\end{Prop}

In the particular case where $r=l$ and $(k_1,\ldots,k_r)=(1,\ldots,l)$, we recover the fact that $G(S)=S$ as shown in \cite[\S 2]{La}. 

\begin{proof}
  Using the hypothesis on the charge and the definition of the action of the Chevalley operators, we see that for all $\underline{i} \in J(S_j)$, if $k\neq j$, we have 
\[
  F_{\underline{i}}\cdot \emptyset_{\mathbf{v}^k}=0.
\]
We thus deduce that for the standard symbol $S=(S_1,\ldots,S_r)$, the element 
\begin{multline*}
  \sum_{(\underline{i}_1,\ldots,\underline{i}_r)\in J(S_1) \times \ldots \times J(S_r)}
    a_{\underline{i}_1}(q)\ldots a_{\underline{i}_r}(q) F_{\underline{i}_1} \ldots F_{\underline{i}_r}\cdot \emptyset_{\mathbf{v}}\\
   \\= \sum_{(T_1,\ldots,T_r)\in \mathcal{S}_{\mathbf{v}^1}(n_1) \times \dots \times \mathcal{S}_{\mathbf{v}^r}(n_r)} \alpha^1_{S_1,T_1}(q)\ldots  \alpha^r_{S_r,T_r}(q) (T_1,\ldots,T_r).
\end{multline*}
is bar invariant and equal to $S$ modulo $q\mathcal{F}_{\mathbf{v},R}$. It is thus the canonical basis element $G(S)$.
\end{proof}

\begin{exa}
  Let us consider $n=4$ and $\mathbf{v} = (6,5,1,0)$. Then $\mathbf{v}^1 = (7,6)$ and $\mathbf{v}^2 = (1,0)$ satisfy the above assumptions. The symbol $S$ given by
  \[
    S=\left(
      \begin{array}{cccccccccc}
        \ldots & -1 & 0 & 1 & 2 & 3 & 4 & 5 & 6 \\
        \ldots & -1 & 0 & 1 & 2 & 3 & 4 & 7 &   \\
        \ldots & -1 & 0 & 2 &   &   &   &   &   \\
        \ldots & -1 & 3 &   &   &   &   &   &   
      \end{array} \right)
  \]
  is standard, of size {$6$} and decomposes as $S=(S_1,S_2)$ with
  \[
    {S_1= \left(
      \begin{array}{ccccc}
        \ldots & 3 & 4 & 5 & 6 \\
        \ldots & 3 & 4 & 7 &   \\
      \end{array} \right)
    \quad\text{and}\quad
    S_2=\left(
      \begin{array}{ccccc}
        \ldots & -2 & -1 & 0 & 2 \\
        \ldots & -2 & -1 & 3 &   \\
      \end{array} \right).}
  \]
  Both $S_1$ and $S_2$ are standard and we have, using for example Leclerc--Miyachi's formula
  \[
    {G(S_1) =
    S_1
    +
    q\left(
      \begin{array}{ccccc}
        \ldots & 3 & 4 & 5 & 7 \\
        \ldots & 3 & 4 & 6 &   \\
      \end{array} \right) 
    \quad\text{and}\quad
    G(S_2) =
    S_2
    +
    q\left(
      \begin{array}{ccccc}
        \ldots & -2 & -1 & 0 & 3 \\
        \ldots & -2 & -1 & 2 &   \\
      \end{array} \right).}
  \]
  Therefore, the canonical basis element $G(S)$ is
  \begin{multline*}
    G(S) =
    S
    +
    q\left(
      \begin{array}{cccccccccc}
        \ldots & -1 & 0 & 1 & 2 & 3 & 4 & 5 & 7 \\
        \ldots & -1 & 0 & 1 & 2 & 3 & 4 & 6 &   \\
        \ldots & -1 & 0 & 2 &   &   &   &   &   \\
        \ldots & -1 & 3 &   &   &   &   &   &   
      \end{array} \right)
    +
    q\left(
      \begin{array}{cccccccccc}
        \ldots & -1 & 0 & 1 & 2 & 3 & 4 & 5 & 6 \\
        \ldots & -1 & 0 & 1 & 2 & 3 & 4 & 7 &   \\
        \ldots & -1 & 0 & 3 &   &   &   &   &   \\
        \ldots & -1 & 2 &   &   &   &   &   &   
      \end{array} \right)
    \\
    +
    q^2\left(
      \begin{array}{cccccccccc}
        \ldots & -1 & 0 & 1 & 2 & 3 & 4 & 5 & 7 \\
        \ldots & -1 & 0 & 1 & 2 & 3 & 4 & 6 &   \\
        \ldots & -1 & 0 & 3 &   &   &   &   &   \\
        \ldots & -1 & 2 &   &   &   &   &   &   
      \end{array} \right).
  \end{multline*}
\end{exa}

\subsection{Column Removal Theorem for symbols}

In this part, we show a "column removal theorem'' for symbols. We fix $\mathbf{v} =(v_1,\ldots,v_l)\in \mathcal{A}^l$ as usual. 

\begin{Th}\label{removal}
  Let $S\in\mathcal{S}^{\st}_{\mathbf{v}}$ be a standard symbol of size $n$. We write the associated canonical basis element as
  \[
    G(S)=\sum_{T \in \mathcal{S}_\mathbf{v} (n) } \alpha_{S,T}(q) T.
  \]
  We assume that there is an integer $x$ that appears in every row of $S$. This common entry is then in every row of every symbol $T$ appearing in $G(S)$ with a {nonzero} coefficient. We denote by $T[x]\in\mathcal{S}_{\mathbf{v}^{-}}$, the symbol obtained from $T$ by removing every entry $x$ from the symbol, where $\mathbf{v}^{-}=(v_1-1,\ldots,v_l-1)$. Then we have 
  \[
    G(S[x])=\sum_{T\in\mathcal{S}_{\mathbf{v}}(n)} \alpha_{S,T}(q) T[x].
  \]
\end{Th}

\begin{proof}
  First, we assume that for all $y<x$ there is an entry $y$ in every row of $S$. Then we write the canonical basis element associated to $S[x]$ as
  \[
    G(S[x]) = \sum_{T}\alpha_{S[x],T}(q)T,
  \]
  the sum being taken on symbols of multicharge $\mathbf{v}^{-}$. Every symbol appearing with a {nonzero} coefficient has the same entries as $S[x]$, and is therefore of the form $T[x]$ for some $T\in \mathcal{S}_{\mathbf{v}}(n)$. We then have
  \[
    G(S[x]) = \sum_{T}\alpha_{S[x],T[x]}(q)T[x],
  \]
  the sum being taken on all the symbols $T\in \mathcal{S}_{\mathbf{v}}(n)$ with the same multiset of entries as $S$. For any $z\in \mathbb{Z}$, let $m_z$ be the number of entries $z$ in the symbol $S$, and let $t$ be the maximal entry. Then, using \Cref{divided}, for all $T[x]$ appearing in the above sum, $E_{x}^{(m_{x+1})}\cdot T[x]$ is the symbol obtained from {$T[x]$} by changing every $x+1$ into $x$, since there is no entry $x$ in $T[x]$. Repeating the same argument, we then have
  \[
    E_{t-1}^{(m_t)}\dots E_{x+1}^{(m_{x+2})}E_x^{(m_{x+1})}\cdot G(S[x]) = \sum_{T}\alpha_{S[x],T[x]}(q)T[x]^{<},
  \]
  where $T[x]^{<}$ is obtained from $T[x]$ by replacing every entry $z>x$ by $z-1$. This element is bar invariant and equal to $S[x]^{<}$ modulo $q\mathcal{F}_{\mathbf{v}^{-},R}${. It is then} the canonical basis element associated to the standard symbol $S[x]^{<}$:
  \[
    G(S[x]^{<}) = \sum_{T}\alpha_{S[x],T[x]}(q)T[x]^{<}.
  \]
  Recall the notation $U^{+}$ of \Cref{plus} for a symbol $U$. Since all $z<x$ is in every row of each symbol $T[x]$ appearing with a {nonzero} coefficient, we have $(T[x]^{<})^{+}=T$. Then
  \[
    G(S) = \sum_{T}\alpha_{S[x],T[x]}(q)T
  \]
  is the expansion in the standard basis of the canonical basis element associated to $S$, which concludes the proof in this case.

  \medskip

  We now proceed by induction on the size of $S$. If $S$ is of size $0$, then $S=\emptyset_{\mathbf{v}}$ and we are in the case treated above.
  
  Let now $S$ be a symbol of size $n$. We may suppose that that there exists $y<x$ and a row of $S$ without the entry $y$. We choose $y$ maximal so that $y+1,y+2,\dots,x$ all appear in every row of $S$. We will first delete the entries $y+1$ and then lower by $1$ every other entries $y+2,\ldots,x$. Let $i_1,\ldots,i_r$ the rows of $S$ containing the entry $y$, with $0\leq r<l$. As before, we write
\[
  G(S) = \sum_{T\in \mathcal{S}_{\mathbf{v}}(n)}\alpha_{S,T}(q)T.
\]
One again, if $\alpha_{S,T}(q)\neq 0$, then the multiset of elements in $S$ and $T$ are the same. Therefore, every $T$ appearing with a {nonzero} coefficient contains the entries $y+1,y+2,\ldots,x$ in every row and the entry $y$ in $r$ rows.
 
For such a symbol $T$, we denote by $\widetilde{T}$ the symbol which is obtained from $T$ by replacing each $y+1$ with $y$ in all the rows where $y$ does not appear. If
\[
  T=\left( 
    \begin{array}{cc}
      (\ast) &  \textcolor{blue}{y+1} \\
      \vdots & \vdots \\
      y &  \textcolor{blue}{y+1} \\
      (\ast)  &  \textcolor{blue}{y+1} \\
      \vdots & \vdots \\
      (\ast) &  \textcolor{blue}{y+1} \\
      y &  \textcolor{blue}{y+1}\\
      \vdots & \vdots \\
      (\ast)  &  \textcolor{blue}{y+1} \\
    \end{array}\right) 
  \begin{array}{c}
    1 \\
    \vdots \\
    i_1\\
    \\
    \vdots \\
    \\
    i_r\\
    \vdots  \\
    l\\
  \end{array}
\]
then
\[
  \widetilde{T}=\left( 
    \begin{array}{cc}
      (\ast) &  \textcolor{blue}{y} \\
      \vdots & \vdots \\
      \textcolor{blue}{y} & y+1 \\
      (\ast) &  \textcolor{blue}{y} \\
      \vdots & \vdots \\
      (\ast) &  \textcolor{blue}{y} \\
      \textcolor{blue}{y}& y+1\\
      \vdots & \vdots \\
      (\ast)  &  \textcolor{blue}{y} \\
    \end{array}\right) 
  \begin{array}{c}
    1 \\
    \vdots \\
    i_1\\
    \\
    \vdots \\
    \\
    i_r\\
    \vdots  \\
    l\\
  \end{array}
\]
In the above symbols, $(\ast)$ denotes an entry smaller than $y$. Moreover, we wrote in the same column the entries $y+1$ in $T$, {even though} they need not be in the same column.

Using \Cref{divided}, we then have  
\[
  E^{(l-r)}_{y}\cdot G(S)=\sum_{T\in \mathcal{S}_{\mathbf{v}}(n)} \alpha_{S,T}(q) E^{(l-r)}_{y}\cdot T =\sum_{T\in \mathcal{S}_{\mathbf{v}}(n)}\alpha_{S,T}(q) \widetilde{T}.
\]
Indeed, that there are exactly $l-r$ elements $y+1$ in the symbol of $\widetilde{T}$ with no $y$ in the same row. From the definition of the canonical basis element associated to the standard symbol $\widetilde{S}$, we deduce that
\[
  G(\widetilde{S}) = E^{(l-r)}_{y}\cdot G(S) =\sum_{T\in\mathcal{S}_\mathbf{v}(n)} \alpha_{S,T}(q) \widetilde{T}.
\]
Now, the symbol $\widetilde{S}$ is of size strictly smaller than $S$ and contains $y$ in every row. By induction, we have
\[
  G(\widetilde{S}[y])=\sum_{T\in \mathcal{S}_{\mathbf{v}}(n)} \alpha_{S,T}(q) \widetilde{T}[y].
\]
Once again, we have 
\[
  E^{(r)}_{y}\cdot G(\widetilde{S}[y])=\sum_{T\in \mathcal{S}_{\mathbf{v}}(n)} \alpha_{S,T} (q) E^{(r)}_{y}\cdot\widetilde{T}[y].
\]
If $\alpha_{S,T}(q) \neq 0$, then there are exactly $r$ element $y+1$ in $\widetilde{T}[y]$ without any $y$ in the symbol. This implies that for all such symbols, we get $E^{(r)}_{y}\cdot{\widetilde{T}[y]}=T[y+1]$ and thus
\[
  E^{(r)}_{y} \cdot G(\widetilde{S}[y])=\sum_{T\in \mathcal{S}_{\mathbf{v}}(n)} \alpha_{S,T} (q) T[y+1].
\]
As this element is bar invariant and equal to $S[y+1]$ modulo $q\mathcal{F}_{\mathbf{v}^{-},R}$, we have
\[
  G(S[y+1]) = E^{(r)}_{y} \cdot G(\widetilde{S}[y]) = \sum_{T\in \mathcal{S}_{\mathbf{v}}(n)} \alpha_{S,T} (q) T[y+1].
\]
We have removed all entries $y+1$ from the symbol $S$ without changing the expression of the canonical basis element. It remains to lower by $1$ every entry $y+2,y+3,\ldots,x$. For all $T[y+1]$ appearing in the above sum, we have $E_{x-1}^{(l)}\dots E_{y+1}^{(l)}E_{y+1}^{(l)}\cdot T[y+1] = T[x]$ and therefore, using the definintion of the canonical basis element, we obtain that
\[
  G(S[x]) = E_{x-1}^{(l)}\dots E_{y+1}^{(l)}E_{y+1}^{(l)}\cdot G(S[y+1]) = \sum_{T\in \mathcal{S}_{\mathbf{v}}(n)} \alpha_{S,T} (q) T[x],
\] 
and the proof is complete.
\end{proof}

\begin{Rem}
  If $l=2$, the above result can be easily deduced from the Leclerc--Miyachi's formula because if $\beta_{j_1}^1=\beta_{j_2}^2$ then $\Psi(\beta_{j_2}^2)=\beta_{j_1}^1$ and the pairs of $S$ and of {$S[x]$} are the same. 
\end{Rem}

\begin{Rem}
{The column removal theorem can be used to compute the canonical basis element $G(S[x])$ from $G(S)$ but also the canonical basis element $G(S)$ from $G(S[x])$. As the calculation of these canonical basis elements are often done by induction on the size of the symbol, it is interesting to notice that the size of $S[x]$ might be smaller or greater than the size of $S$.}
\end{Rem}

{The examples below are computed using the algorithm developed in \cite{Jal}.}

  \begin{exa}\label{l1}
Assume that $l=3$, $\mathbf{v}=(2,2,1)$ and take the standard symbol:
\[
  S=\left( 
    \begin{array}{ccccc}
      \ldots & 0 & 2 & 3 \\
      \ldots & 0 & 2 & 4 \\
      \ldots & 0 & 2  
    \end{array}\right).
\]
Removing all occurences of the entry $2$ gives the symbol
\[
  S[2]=\left( 
    \begin{array}{ccccc}
      \ldots & 0 & 3  \\
      \ldots & 0 & 4  \\
      \ldots & 0
    \end{array}\right).
\]
Note that the associated multicharge is now $(1,1,0)$. Now we have:
\[
  G(S[2])=\left( 
    \begin{array}{cccc}
      \ldots & 0 & 3 \\
      \ldots & 0 & 4 \\
      \ldots & 0 &  
    \end{array}\right)
  +q
  \left( 
    \begin{array}{ccccc}
      \ldots & 0 & 4 \\
      \ldots & 0 & 3 \\
      \ldots & 0  
    \end{array}\right)
\]
and 
\[
  G(S)=\left( 
    \begin{array}{ccccc}
      \ldots & 0 & 2 & 3 \\
      \ldots & 0 & 2 & 4 \\
      \ldots & 0 & 2
    \end{array}\right)
  +q
  \left( 
    \begin{array}{ccccc}
      \ldots & 0 & 2 & 4 \\
      \ldots & 0 & 2 & 3 \\
      \ldots & 0 & 2 
    \end{array}\right)
\]
which is consistent with our result. {In this example, the size of $S[2]$ is $5$ and the size of $S$ is $6$.}
\end{exa}

\begin{exa}
  Assume that $l=4$, $\mathbf{v}=(4,3,1,1)$ and take the standard symbol:
  \[
    S=\left( 
      \begin{array}{cccccc}
        \ldots & 0 & 1 & 2 & 3 & 5 \\
        \ldots & 0 & 1 & 3 \\
        \ldots & 0 & 1 \\
        \ldots & 1 & 2
      \end{array}\right).
  \]
  Removing all occurences of the entry $1$ gives the symbol
  \[
    S[1]=\left( 
      \begin{array}{cccccc}
        \ldots & 0 & 2 & 3 & 5 \\
        \ldots & 0 & 3 \\
        \ldots & 0 \\
        \ldots & 2
      \end{array}\right).
  \]
  Note that the associated multicharge is now $(3,1,0,0)$ Now we have:
  \[
    G(S[1])=\left(
      \begin{array}{cccccc}
        \ldots & 0 & 2 & 3 & 5 \\
        \ldots & 0 & 3 \\
        \ldots & 0 \\
        \ldots & 2
      \end{array}\right)
    +q
    \left( 
      \begin{array}{cccccc}
        \ldots & 0 & 2 & 3 & 5 \\
        \ldots & 0 & 3 \\
        \ldots & 2 \\
        \ldots & 0
      \end{array}\right)
    +q^2\left(
      \begin{array}{cccccc}
        \ldots & 0 & 2 & 3 & 5 \\
        \ldots & 2 & 3 \\
        \ldots & 0 \\
        \ldots & 0
      \end{array}\right)
  \]
  and 
  \[
    G(S)=\left(
      \begin{array}{cccccc}
        \ldots & 0 & 1 & 2 & 3 & 5 \\
        \ldots & 0 & 1 & 3 \\
        \ldots & 0 & 1 \\
        \ldots & 1 & 2
      \end{array}\right)
    +q
    \left( 
      \begin{array}{cccccc}
        \ldots & 0 & 1 & 2 & 3 & 5 \\
        \ldots & 0 & 1 & 3 \\
        \ldots & 1 & 2 \\
        \ldots & 0 & 1
      \end{array}\right)
    +q^2
    \left(
      \begin{array}{cccccc}
        \ldots & 0 & 1 & 2 & 3 & 5 \\
        \ldots & 1 & 2 & 3 \\
        \ldots & 0 & 1 \\
        \ldots & 0 & 1
      \end{array}\right)
  \]
  which is consistent with our result. {In this example, the size of $S[1]$ is $8$ and the size of $S$ is $4$.}
\end{exa}

\section{Computation of  canonical basis elements}

We investigate two particular cases in which the elements of the canonical basis can be explicitly computed. In each of these cases, we provide closed formulae for the canonical basis elements, which also prove that they are monomial, in the same spirit as the result of Leclerc--Miyachi. 

\subsection{Ordered symbol}

Assume as usual that $\mathbf{v}\in \mathcal{A}^l$. We now define a certain distinguished class of symbol. We say that a symbol $S=(\beta^1,\ldots,\beta^l)$ written as follows
\[
  S =
  \begin{pmatrix}
    \beta^1 \\
    \beta^{2} \\
    \vdots \\
    \beta^l
  \end{pmatrix}
  =
  \begin{pmatrix}
    \cdots & \beta^{1}_{v_l-1} & \beta^{1}_{v_l} & \beta^{1}_{v_l+1} & \cdots & \beta^{1}_{v_{2}-1} &  \beta^{1}_{v_{2}} & \beta^{1}_{v_{2}+1} & \cdots & \beta^{1}_{v_1} \\
    \cdots & \beta^{2}_{v_l-1} & \beta^{2}_{v_l} & \beta^{2}_{v_l+1} & \cdots & \beta^{2}_{v_{2}-1} &  \beta^{2}_{v_{2}} &  &  &  \\
    \vdots & \vdots & \vdots & \vdots & \vdots & \vdots & \vdots & & & \\
    \cdots & \beta^l_{v_l-1} & \beta^l_{v_l} & & & & & & & 
  \end{pmatrix}
\]
is an {\it ordered symbol} if for all relevant $c \in \{1,\ldots,l-1\}$ and  relevant $k\in \mathbb{Z}$, we have $\beta^c_k \leq \beta^{c+1}_{k}$  (that is, $S$ is standard) and if for all $i\leq v_{1}$ and for all relevant $c \in \{1,\ldots,l\}$, we have that $\beta^{1}_i \geq \beta^c_{i-1}$. In other words, the entries of the symbol {increase} from top to bottom and left to right. Note that in each block there exists a unique such symbol. 

\begin{exa}
  The following symbol is an ordered symbol associated with the multicharge $(2,2,1)$:
  \[
    S=\left( 
      \begin{array}{ccccc}
        \ldots & -1 & 0 & 1 & 5 \\
        \ldots & -1 & 0 & 2 & 5 \\
        \ldots & -1 & 1 & 4
      \end{array}
    \right).
  \]
\end{exa}
 
For $r\in \mathbb{N}$, we consider the symmetric group $\mathfrak{S}_r$, we write $s_j$ for the transposition $(j,j+1)$ with $j=1,\ldots,r-1$ and denote by $\ell$ the length function with respect to these generators. Let $S=(\beta^1,\ldots,\beta^l)$ be an ordered symbol with the same convention as above. Define $h_j$ to be the height of the $j$th column of $S$. Then $h_j=r$ if and only if {$v_{r+1}< j \leq v_r$}. Let us denote {the $j$th column} of the symbol as follows:
\[
  X_j=(\beta_j^{1},\ldots, \beta_j^{h_j}).
\]
For $\sigma \in \mathfrak{S}_{h_j}$, we denote $(X_j)^{\sigma}=(\beta_j^{\sigma (1)},\ldots, \beta_j^{\sigma(h_j)})$ the permuted column.

Let $\mathfrak{S}_{\mathbf{v}}=\prod_{j\leq v_{1}}\mathfrak{S}_{h_j}$. The group $\mathfrak{S}_{\mathbf{v}}$ acts on symbols in $\mathcal{S}_{\mathbf{v}}$ by permuting the entries columnwise. The length $\ell(\sigma) \in \mathbb{N}\cup\{+\infty\}$ of a tuple $\sigma=(\sigma_{j})_{j\leq v_{1}}$ is the sum of the length of its components $\sigma_j$. Here, we allow $\ell(\sigma)$ to be infinite, but will not encounter elements of $\mathfrak{S}_{\mathbf{v}}$ with infinite length.

We say that $\sigma$ is admissible for the symbol $S\in \mathcal{S}_{\mathbf{v}}$ if 
\begin{enumerate}[label=\roman*)]
\item The following 
  \[
    \begin{pmatrix}
      \cdots & \beta^{\sigma_{v_l-1}(1)}_{v_l-1} & \beta^{\sigma_{v_l}(1)}_{v_l} & \beta^{\sigma_{v_l+1}(1)}_{v_l+1} & \cdots & \beta^{\sigma_{v_{2}-1}(1)}_{v_{2}-1} &  \beta^{\sigma_{v_{2}}(1)}_{v_{2}} & \beta^{\sigma_{v_{2}+1}(1)}_{v_{2}+1} & \cdots & \beta^{\sigma_{v_{1}(1)}}_{v_1} \\
      \cdots & \beta^{\sigma_{v_l-1}(2)}_{v_l-1} & \beta^{\sigma_{v_l}(2)}_{v_l} & \beta^{\sigma_{v_l+1}(2)}_{v_l+1} & \cdots & \beta^{\sigma_{v_{2}-1}(2)}_{v_{2}-1} &  \beta^{\sigma_{v_{2}}(2)}_{v_{2}} &  &  &  \\
      \vdots & \vdots & \vdots & \vdots & \vdots & \vdots & \vdots & & & \\
      \cdots & \beta^{\sigma_{v_l-1}(l)}_{v_l-1} & \beta^{\sigma_{v_l}(l)}_{v_l} & & & & & & & 
    \end{pmatrix}
  \]
  is a well-defined symbol, that is there are no repetition in each row. We denote this symbol by $S^{\sigma}$. 
  
\item $\sigma$ has minimal length among all the $w \in \mathfrak{S}_{\mathbf{v}}$ such that  $S^{\sigma}=S^{w}$.
\end{enumerate}

Note that, since for every symbol the $j$-th column is eventually $(j,\ldots,j)$, the second condition ensures that an admissible permutation $\sigma$ is of finite length, since $\sigma_j=\id$ for $j$ small enough.
 
\begin{exa}
  We keep the symbol $S$ of the previous example. Then $\sigma=(\dots,\id,s_2,s_2,\id)\in \mathfrak{S}_{\mathbf{v}}$ is admissible for $S$, and we have 
  \[
    S^{\sigma}=\left( 
      \begin{array}{ccccc}
        \ldots & -1 & 0 & 1 & 5\\
        \ldots & -1 & 1 & 4 & 5\\
        \ldots & -1 & 0 & 2
      \end{array}
    \right)
  \]
  we see that this is indeed a well-defined symbol. If now we take $\sigma=(\dots,\id,s_2,s_1,\id)\in\mathfrak{S}_{\mathbf{v}}$ then 
  \[
    S^{\sigma}=\left( 
      \begin{array}{ccccc}
        \ldots & -1 & 0 & 2 & 5\\
        \ldots & -1 & 1 & 1 & 5\\
        \ldots & -1 & 0 & 4
      \end{array}
    \right)
  \]
  is not a well-defined symbol because we have a repetition in the second row. So $\sigma$ is not admissible.
\end{exa}

We denote by $\mathfrak{S}(S)$ the set of elements in $\mathfrak{S}_{\mathbf{v}}$ admissible for the symbol $S$. 
Now for a symbol $T=(\beta^1,\ldots,\beta^l)$,  we define the following number:
\[
  M(T) =\sum_{j \leq v_l} \sum_{t=1}^{h_j} \sharp \left\{k\in \{1,\ldots,t-1\} \middle\vert\ \beta_{j-1}^{k} =\beta_j^{t} \right\}{.}
\]
Note that only there is a finite number of {nonzero} terms in the above sum. 

\begin{Th}\label{ordered}
  Let $S$ be an ordered symbol. Then the canonical basis element is given by
  \[
    G(S)= \sum_{\sigma \in \mathfrak{S}(S)} q^{N(\sigma,S)} S^{\sigma},
  \]
  where $N(\sigma,S)=\ell(\sigma)-M(S^{\sigma})$. In addition, this canonical basis element is monomial.
\end{Th}

\begin{exa}
  Let $\mathbf{v}=(2,2,1)$ and take the standard symbol:
  \[
    S= \left(
      \begin{array}{ccccc}
        \ldots & 0 & 1 & 3 & 5 \\
        \ldots & 0 & 2 & 3 & 5 \\
        \ldots & 1 & 3 & 4
      \end{array} \right)
  \]
  Then the canonical basis element $G(S)$ is equal to
  \begin{multline*}
    \left(
      \begin{array}{ccccc}
        \ldots & 0 & 1 & 3 & 5 \\
        \ldots & 0 & 2 & 3 & 5 \\
        \ldots & 1 & 3 & 4
      \end{array} \right)
    +q
    \left(
      \begin{array}{ccccc}
        \ldots & 0 & 1 & 3 & 5 \\
        \ldots & 0 & 3 & 4 & 5 \\
        \ldots & 1 & 2 & 3
      \end{array} \right) 
    +q
    \left(
      \begin{array}{ccccc}
        \ldots & 0 & 1 & 3 & 5 \\
        \ldots & 1 & 2 & 3 & 5 \\
        \ldots & 0 & 3 & 4
      \end{array} \right)
    +q
    \left(
      \begin{array}{ccccc}
        \ldots & 0 & 1 & 3 & 5 \\
        \ldots & 1 & 3 & 4 & 5 \\
        \ldots & 0 & 2 & 3
      \end{array} \right)
    \\
    +q
    \left(
      \begin{array}{ccccc}
        \ldots & 0 & 2 & 3 & 5 \\
        \ldots & 0 & 1 & 3 & 5 \\
        \ldots & 1 & 3 & 4
      \end{array} \right)
    +q^2
    \left(
      \begin{array}{ccccc}
        \ldots & 0 & 3 & 4 & 5 \\
        \ldots & 0 & 1 & 3 & 5 \\
        \ldots & 1 & 2 & 3
      \end{array} \right)
    +q^2
    \left(
      \begin{array}{ccccc}
        \ldots & 1 & 2 & 3 & 5 \\
        \ldots & 0 & 1 & 3 & 5 \\
        \ldots & 0 & 3 & 4
      \end{array} \right)
    +q^3
    \left(
      \begin{array}{ccccc}
        \ldots & 1 & 3 & 4 & 5 \\
        \ldots & 0 & 1 & 3 & 5 \\
        \ldots & 0 & 2 & 3
      \end{array} \right)
    \\
    +q^2
    \left(
      \begin{array}{ccccc}
        \ldots & 0 & 2 & 3 & 5 \\
        \ldots & 1 & 3 & 4 & 5 \\
        \ldots & 0 & 1 & 3
      \end{array} \right)
    +q^3
    \left(
      \begin{array}{ccccc}
        \ldots & 0 & 3 & 4 & 5 \\
        \ldots & 1 & 2 & 3 & 5 \\
        \ldots & 0 & 1 & 3
      \end{array} \right)
    +q^3
    \left(
      \begin{array}{ccccc}
        \ldots & 1 & 2 & 3 & 5 \\
        \ldots & 0 & 3 & 4 & 5 \\
        \ldots & 0 & 1 & 3
      \end{array} \right)
    +q^4
    \left(
      \begin{array}{ccccc}
        \ldots & 1 & 3 & 4 & 5 \\
        \ldots & 0 & 2 & 3 & 5 \\
        \ldots & 0 & 1 & 3
      \end{array} \right),
\end{multline*}
where below, we give the triplet $(\sigma, \ell(\sigma),M(S^{\sigma}))$ corresponding to each term, and omit the components of $\sigma_j$ for $j<0$:
\begin{gather*}
  ((\id,\id,\id,\id) ,0,0), ((\id,s_2,s_2,\id),2,1), ((s_2,\id,\id,\id),1,0), ((s_2,s_2,\id,\id),2,1),\\
  ((\id,s_1,\id,\id),1,0), ((\id,s_1s_2,s_1s_2,\id),4,2), ((s_1s_2,s_1,\id,\id),3,1), ((s_1s_2,s_1s_2,s_1s_2,\id),6,3),\\
  ((s_2,s_2s_1,s_2,\id), 4,2),((s_2,s_2 s_1 s_2, s_1 s_2,\id),6,3),((s_1 s_2,s_2 s_1,s_1 s_2,\id),6,3),((s_1 s_2,s_1s_2 s_1,s_1s_2,\id),7,3))
\end{gather*}
\end{exa}

Our \Cref{ordered} will be a {consequence} of a few intermedate results. The strategy is to show that the formula for $G(S)$ is given by a monomial element, that we define inductively.

\begin{lemma}
  Assume $S$ is a standard symbol of {nonzero} size. Then there exists at least one element among 
  \[
    \{\ldots, \beta^l_{v_l-1},\beta^l_{v_l}, \beta^{l-1}_{v_{l}+1},\ldots \beta^{l-1}_{v_{l-1}},\ldots, \beta^{1}_{v_{2}+1},\ldots \beta^{1}_{v_{1}}\},
  \]
  that we denote by $x=\beta^k_i$, such that $x-1\neq \beta^s_j$ for all $s\geq k$ and $j<i$.
\end{lemma}

\begin{proof}
Consider the symbol 
\[
  S=\begin{pmatrix}
    \cdots & \beta^{1}_{v_l-1} & \beta^{1}_{v_l} & \beta^{1}_{v_l+1} & \cdots & \beta^{1}_{v_{2}-1} & \beta^{1}_{v_{2}} & \textcolor{red}{\beta^{1}_{v_{2}+1}} & \cdots & \textcolor{red}{\beta^{1}_{v_1}} \\
    \cdots & \beta^{2}_{v_l-1} & \beta^{2}_{v_l} & \beta^{2}_{v_l+1} & \cdots & \beta^{2}_{v_{2}-1} & \textcolor{red}{ \beta^{2}_{v_{2}} }&  &  &  \\
    \vdots & \vdots & \vdots & \vdots & \vdots & \vdots & \vdots & & & \\
    \cdots & \textcolor{red}{\beta^l_{v_l-1} }& \textcolor{red}{\beta^l_{v_l}}& & & & & & & 
  \end{pmatrix}.
\]

We start with the set of numbers $(\beta^l_j)_{j\leq v_l}$. If there exists an element $\beta^l_j$ such that $\beta_{j-1}^l\neq \beta^l_j-1$ then we can set {$x=\beta^l_j$} and we are done. Otherwise, we switch to the second subset of $\beta$-numbers. If among $(\beta^{l-1}_j)_{v_l<j\leq v_{l-1}}$ there is an entry as in the conclusion of the lemma, we are done. Otherwise, we show that $\beta_{j}^{l-1}=j$ for all $j\leq v_{l-1}$. Since $S$ is standard, we have, for all $j\leq v_{l}$, $\beta_{j}^{l-1}\leq \beta_{j}^l$. But we always have for all $1\leq k\leq l$ and $j\leq v_{k}$, $\beta_{j}^k\geq j$. Therefore, $\beta_{j}^{l-1}=j$ for all $j\leq v_{l}$. Now, this also forces $\beta_j^{l-1}=j$ for $v_{l}<j\leq v_{l-1}$, since there are no element among $(\beta^{l-1}_j)_{v_l< j\leq v_{l-1}}$ satisfyting the conclusion of the lemma. Continuing in this way, the process ends otherwise, the symbol $S$ is the symbol
\[
  S=\begin{pmatrix}
    \cdots & v_l-1 & v_l & v_l+1 & \cdots & v_{2}-1 & v_{2} & v_{2}+1 & \cdots & v_1 \\
    \cdots & v_l-1 & v_l & v_l+1 & \cdots & v_{2}-1 & v_{2}&  &  &  \\
    \vdots & \vdots & \vdots & \vdots & \vdots & \vdots & \vdots & & & \\
    \cdots & v_l-1& v_l& & & & & & & 
  \end{pmatrix}
\]
which is $\emptyset_{\mathbf{v}}$ of size $0$.
\end{proof}

Given an ordered symbol $S$, we define a sequence $\seq(S)$ of integers as follows. Let $x=\beta_i^k$ given by above lemma, chosen such that first is $k$ is maximal and then $i$ is minimal among the valid choices. Assume in addition that
\[
  x=\beta_i^k=\beta_i^{k-1}=\ldots =\beta_i^{k-a+1}
\]
with $a$ maximal. Then, the only entries in the symbol $S$ which could be equal to $x-1$ are of the form $\beta_i^s$ for $s<k+a-1$.
 
If we replace by $x-1$ all the elements $\beta_i^{k}=\ldots =\beta_i^{k+a-1}=x$ in $S$, we obtain a well defined symbol $S'$ which is clearly ordered. We then define by induction the sequence $\seq(S)$ by
\[
  \seq (S)=\seq (S'),\underbrace{x\ldots,x}_{\text{ $a$ times}}.
\]

\begin{exa}\label{l1}
  Take the following ordered symbol:
  \[
    \left(
      \begin{array}{ccccccc}
        \ldots & 0 & 1 & 2 & 4 & 5 & 7\\
        \ldots & 0 & 2 & 3 & 5 \\
        \ldots & 0 & 2 & 4 \\
        \ldots & 0 & 2 & 4
      \end{array}\right).
  \]
  Then, we take $x=2$ and the last three elements of the sequence $\seq(S)$ are $2,2,2$. We obtain the symbol
  \[
    \left(
      \begin{array}{ccccccc}
        \ldots & 0 & 1 & 2 & 4 & 5 & 7\\
        \ldots & 0 & 1 & 3 & 5 \\
        \ldots & 0 & 1 & 4\\
        \ldots & 0 & 1 & 4
      \end{array}\right),
  \]
  and now take $x=4$. Then the sequence $\seq(S)$ ends with $4,4,2,2,2$. We repeat this process until we obtain the symbol $\emptyset_{\mathbf{v}}$ of size $0$, and the sequence $\seq(S)$ is {$6,7,5,4,4,5,3,3,3,4,4,2,2,2$}.
\end{exa}
  
Now assume that $\seq(S)$ is
\[
  \underbrace{x_1,\ldots,x_1}_{\text{$a_1$ times}},\underbrace{x_2,\ldots,x_2}_{\text{$a_2$ times}},\ldots,\underbrace{x_m,\ldots,x_m}_{\text{$a_m$ times}}.
\]
We then set
\[
  A(S)=F^{(a_m)}_{x_m-1} \ldots F^{(a_1)}_{x_1-1}\cdot \emptyset_{\mathbf{v}}.
\]

\begin{Prop}\label{ordered_monomial}
  Let $S$ be an ordered symbol. Then we have 
  \[
    A(S)= \sum_{\sigma\in\mathfrak{S}(S)} q^{N(\sigma,S)} S^{\sigma}
  \]
  where $N(\sigma,S)=\ell(\sigma)-M(S^{\sigma})$.
\end{Prop}

\begin{proof}
  This is done by induction on $m$. The result is trivial for $m=0$. Assume that, keeping the above notations, we have 
  \[
    A(S')= \sum_{\sigma\in \mathfrak{S}(S')} q^{N(\sigma,S')}(S')^{\sigma}.
  \]
  We will write $x$ instead of $x_m$. Consider the column $X_j$ of $S$ which contains $x$ on its last line.  We will write $h$ instead of $h_j$ for the height of $X_j$.

One can assume that there are $a_m$ coefficient $x$ on the last $a_m$ rows of the column $X_j$ and then $b$ coefficient $x-1$ with $b\geq 0$ on the same column $X_j$. We will write the colum $X_j$ as follows:
\[
  X_j={(y_1,\ldots,y_{h-b-a_m}, \underbrace{x-1,\ldots,x-1}_{b\text{ times}} ,\underbrace{x,\ldots,x}_{a_m\text{ times}})},
\]
where we denoted by $h$ the height of the column and by ${y_1,\ldots,y_{h-b-a_m}}$ the entries of $X_j$ smaller than $x-1$.

In the column $X_{j-1}$, this is not possible to have any $x$, by the choice of $x$, nor $x-1$, since $S$ is ordered. In the column $X_{j+1}$, we may have several elements equal to $x$, at the top of the column. The symbol $S'$ has all its columns equal to the columns of $S$ except for the $j$-th column $X'_{j}$, which is then
\[
  X'_{j}={(y_1,\ldots,y_{h-b-a_m}, \underbrace{x-1,\ldots,x-1}_{b\text{ times}} ,\underbrace{x-1,\ldots,x-1}_{a_m\text{ times}}).}
\]

Let $\sigma=(\sigma_{k})_{k} \in \mathfrak{S}(S')$. In $A(S')$, the coefficient of $(S')^{\sigma}$ is $q^{N(\sigma,S')}$. In $X'_{j}$, we have exactly $a_m+b$ coefficients $x-1$. After applying $F^{(a_m)}_{x-1}$, only the entries in column $j$ will change because this is the only column with an entry $x-1$. 
 
We now compute $F^{(a_m)}_{x-1}\cdot(S')^{\sigma}$. The symbols appearing in the decomposition are obtained by changing $a_m$ coefficients $x-1$ into $x$. Let us fix one of these symbols $T$. To obtain explicitly this symbol, we have to consider the sets $Y_k=(X'_k)^{\sigma_k}$ for $k\neq j$ and the set $Y_j$ which is obtained from $(X'_j)^{\sigma_{j}}$ by replacing $a_m$ elements $x-1$ with $x$ (in rows $j_1<\ldots<j_{a_m}$), so that the resulting set have exactly $a_m$ entries $x$ and $b$ entries $x-1$. Let $T$ be the resulting symbol.
 
This coefficient of $T$ in the vector $F^{(a_m)}_{x-1}\cdot(S')^{\sigma}$ is $q^{N(\sigma,S')+N}$, where by \Cref{divided}, $N=\sum_{1\leq i \leq a_m} N_i$ with $N_i$ is the number $N_i^{-}$ of $x-1$ (in $Y_j$) in row $>j_i$ minus the number $N_i^+$ of $x$ (in $Y'_j$) in row $>j_i$. 

For $k\in\{1,\ldots,r-1\}$, let us denote 
\[
  W_k=\{1,s_k,s_k s_{k-1},\ldots,s_k s_{k-1} \ldots s_{1}\},
\]
where $s_j$ corresponds to the transposition $(j,j+1)$. Recall that each element of {$\sigma\in\mathfrak{S}_r$} can be uniquely written as {$w_1w_2\dots w_{r-1}$ where each $w_i$ is an element of $W_i$, and such an expression yields a reduced form for $\sigma$.}
 
Let us consider now $X_j$ and act by the following element on this column:
\[
  \tau= \underbrace{s_{h-a_m}s_{h-a_m-1}\ldots s_{h-N_{1}^{-}-a_{m}+1}}_{\in W_{h-a_m}}\underbrace{s_{h-a_m-1}s_{h-a_m-2}\ldots s_{h-N_{2}^{-}-a_m+2}}_{\in {W_{h-a_{m}-1}}}\ldots\underbrace{{s_{h-1}s_{h-2}}\ldots s_{h-N_{a_m}^{-}}}_{\in {W_{h-1}}}.
\]
This is an element in $W_{h-a_m}\ldots W_{h-1}W_{h}$ which is reduced and of length $\sum_{1\leq i\leq a_m} N_i^-$. We have 
\[
  (X_j)^{\tau}={(y_1,\ldots,y_{h-b-a_m},\underbrace{x-1,\ldots,x-1}_{b-N_1^-},x,\underbrace{x-1\ldots,x-1}_{N_1^{-}-N_{2}^{-}},x,x-1\ldots,x-1,x,\underbrace{x-1\ldots,x-1}_{N^-_{a_m}})}.
\]
We still denote by $\tau$ the element of $\mathfrak{S}_{\mathbf{v}}$ with entry $j$ equal to $\tau$ and all others equal to $\id$. It is then easy to see that $\tau\sigma\in \mathfrak{S}(S)$ and that $S^{\tau\sigma}=T$. Moreover we have 
\[
  \ell(\tau\sigma)=\ell(\sigma)+\sum_{1\leq i\leq a_m} N_i^{-}.
\]
Now the number $\sum_{1\leq i \leq a_m} N_i^+$ corresponds to the number  $M(S^{\tau\sigma})-M({S'}^{\sigma})$ (this is the contribution of elements of the form $x$ in column $j+1$). To summarize, the coefficient of $T$ is $q^{\ell(\sigma)-M({S'}^{\sigma})+N}$ and
\begin{align*}
  \ell(\sigma)-M({S'}^{\sigma})+N
  &= \ell(\sigma)-M({S'}^{\sigma})+\sum_{1\leq i\leq a_m} (N_i^{-}-N_i^{+})\\
  &= \ell(\sigma)+\sum_{1\leq i \leq a_m}N_i^{-}-M({S'}^{\sigma}) - \sum_{1\leq i \leq a_m}N_i^{+}\\
  &= \ell(\tau\sigma) -M(S^{\tau\sigma})
\end{align*}
which is exactly what we wanted. 

In addition, it is clear that any element of the form $S^{\sigma}$ can be obtained from an element ${S'}^{\sigma'}$ by replacing $a_m$ elements $x-1$ by $x$. The result thus follows. 
\end{proof}

\begin{Prop}\label{ordered_positive}
  Assume that $S$ is an ordered symbol and that $\sigma\in \mathfrak{S}(S)$. Then we have $N(\id,S)=0$ and $N(\sigma,S)>0$ if $\sigma \neq \id$.
\end{Prop}

\begin{proof}
The fact that $N(\id,S)=0$ directly follows from the definition of the ordered symbols because $M(S)=0$. 

Assume now that $\sigma=(\sigma_j)_{j\leq v_l}\neq\id$. We claim that, for every $j\leq v_1$, we have
\begin{multline*}
  \ell(\sigma_j) \geq \sum_{t=1}^{h_{j}}\sharp\left\{k\in \{1,\ldots,t-1\} \middle\vert\ \beta_{j-1}^{\sigma_{j-1}(k)} =\beta_j^{\sigma_j(t)} \right\} \quad\text{and}\\ \ell(\sigma_{j-1}) \geq \sum_{t=1}^{h_{j}}\sharp \left\{k\in \{1,\ldots,t-1\} \middle\vert\ \beta_{j-1}^{\sigma_{j-1}(k)} =\beta_j^{\sigma_j(t)} \right\}.
\end{multline*}

Since the symbol $S$ is ordered, each pair of consecutive columns $X_{j-1}$ and $X_{j}$ of $S$ have at most one common value, which would be at the bottom of $X_{j-1}$ and at the top of $X_{j}$. If there is no common value, then the claim is trivial. Otherwise, denote by $x$ the common element of $X_{j-1}$ and of $X_j$, which is then the top entry of the $j$th column $X_{j}$ and the bottom entry of $X_{j-1}$, and suppose that this entry appears $r$ times in $X_j$. Denote by $1\leq i_1<\ldots<i_r \leq h_j$ the rows in which there is an entry $x$ in $X_j^{\sigma_j}$. Then
\[
  \sum_{t=1}^{h_{j}}\sharp\left\{k\in \{1,\ldots,t-1\} \middle\vert\ \beta_{j-1}^{\sigma_{j-1}(k)} =\beta_j^{\sigma_j(t)} \right\} = \sum_{s=1}^{r}\sharp\left\{k\in \{1,\ldots,i_s-1\} \middle\vert\ \beta_{j-1}^{\sigma_{j-1}(k)} = x \right\}.
\]
To move the entries $x$ in rows $1<\ldots<r$ in $X_j$ to the rows $i_1<\ldots<i_r$ in $X_j^{\sigma_j}$, the permutation $\sigma_j$ must be at least of length
\[
  (i_r-r)+(i_{r-1}-(r-1))+\dots+(i_1-1) = \sum_{s=1}^{r}(i_s-s).
\]
Then, since $\sigma$ is admissible, there is no $x$ in $X_{j-1}^{\sigma_{j-1}}$ in the rows $i_1,\dots,i_r$ and we have
\[
  \sharp \left\{k\in \{1,\ldots,i_s-1\} \middle\vert\ \beta_{j-1}^{\sigma_{j-1}(k)} = x \right\} \leq i_s - s,
\]
because among the $i_s$ first rows, there could be $x$ in any rows except $i_1,\ldots,i_s$. Therefore
\begin{align*}
  \sum_{t=1}^{h_{j}}\sharp\left\{k\in \{1,\ldots,t-1\} \middle\vert\ \beta_{j-1}^{\sigma_{j-1}(k)} =\beta_j^{\sigma_j(t)} \right\}
  &= \sum_{s=1}^{r}\sharp\left\{k\in \{1,\ldots,i_s-1\} \middle\vert\ \beta_{j-1}^{\sigma_{j-1}(k)} = x \right\}\\
  &\leq \sum_{s=1}^r (i_s-s)\\
  &\leq \ell(\sigma_j).
\end{align*}
The other inequality is proven similarly.

Finally, if we choose $j\leq v_1$ minimal such that $\sigma_{j-1}\neq\id$, we obtain
\[
  \ell(\sigma_j)+\dots+\ell(\sigma_{v_1}) \geq {M(S^\sigma)},
\]
and therefore $N(S,\sigma)\geq \ell(\sigma_{j-1})>0$.
\end{proof}

\begin{proof}[Proof of \Cref{ordered}]
  From \Cref{ordered_monomial}, the element
  \[
    A(S) = \sum_{\sigma\in\mathfrak{S}(S)} q^{N(\sigma,S)} S^{\sigma}
  \]
  is bar invariant, since it is monomial. Moreover, \Cref{ordered_positive} shows that $A(S) = S \mod q\mathcal{F}_{\mathbf{v},R}$. Then $A(S)$ is the canonical basis element $G(S)$ indexed by the standard symbol $S$. In particular, $G(S)$ is monomial.
\end{proof}

\subsection{A generalization of Leclerc--Miyachi's formula}

We now turn to a last formula which may be viewed as a generalization of the Leclerc--Myiachi's one. Let $S=(\beta^1,\ldots,\beta^l)$ be a standard $l$-symbol of multicharge $\mathbf{v}$. For any $1 \leq i < j \leq l$, let us consider the  map ${\psi_{j,i}} \colon \beta^j \rightarrow \beta^i$ defined as in \eqref{Lmy}.

We consider the following condition \eqref{eq:condition} on standard $l$-symbols:
\begin{equation}
  \label{eq:condition}
  \forall 1 \leq i < j < k \leq l, \ \psi_{k,i} = \psi_{j,i}\circ \psi_{k,j}. \tag{$\heartsuit$}
\end{equation}

\begin{exa}
  Let $\mathbf{v}=(3,3,2)$ and consider the symbols
  \[
    S_1=\left(
      \begin{array}{ccccc}
        \ldots & 0 & 1 & 3 & 5 \\
        \ldots & 0 & 2 & 3 & 5 \\
        \ldots & 1 & 3 & 4
      \end{array}
    \right)
    ,\quad
    S_2=\left(
      \begin{array}{ccccc}
        \ldots & 0 & 2 & 3 & 5 \\
        \ldots & 0 & 2 & 3 & 5 \\
        \ldots & 1 & 3 & 4
      \end{array}
    \right)
    \quad\text{and}\quad
    S_3=\left(
      \begin{array}{ccccc}
        \ldots & 0 & 1 & 3 & 5 \\
        \ldots & 0 & 1 & 3 & 5 \\
        \ldots & 2 & 3 & 4
      \end{array}
    \right).
  \]
  Then the above condition \eqref{eq:condition} is not satisfied by $S_1$ since we have $\psi_{3,2}(1)=0$, $\psi_{2,1}(0)=0$ but $\psi_{3,1}(1)=1$. One may check that the condition \eqref{eq:condition} is satisfied by $S_2$ and $S_3$.
\end{exa}

\begin{Rem}\label{ord}
  Assume that $S$ is an ordered symbol and that for all columns of the symbol, the last element (that is the element at the bottom of it) is different from the first element (that is the element at the top) of the next column. Then the condition \eqref{eq:condition} is satisfied since for every $i\leq j$ and $k\leq v_j$ we have $\psi_{j,i}(\beta_{k}^j) = \beta_{k}^{i}$.

  We have given above an example of an ordered symbol which does not satisfied \eqref{eq:condition}, an example of an ordered symbol satisfying \eqref{eq:condition} but not the extra assumption of the remark, and a non ordered symbol satisfying the condition \eqref{eq:condition}.
\end{Rem}

Note that the condition \eqref{eq:condition} is automatically satisfied in the level $2$ case. The aim of this section is to show that one can obtain a formula for the canonical basis which is a generalization of Leclerc--Miyachi's formula for $l=2$.

We define an equivalence relation on the entries of a symbol $S=(\beta^1,\ldots,\beta^l)$ satisfying the condition \eqref{eq:condition}. We say that $x \in \beta^i$ and $y \in \beta^j$ are in relation with $i<j$ if $\psi_{j,i} (y)=x$. The equivalence classes under this relation are called the {\it spines} of the symbol. Since $S$ satisfy the condition \eqref{eq:condition}, each spine contains at most one element from each row. Abusing notation, we denote them as a tuple of integers starting with the element in $\beta^1$ ending with $\beta^{k}$ with $k$ maximal, and we say that $k$ is the length of such a spine. Then a spine of length $k$ is of the form $(\psi_{k,1}(x),\psi_{k,2}(x),\ldots,\psi_{k,k-1}(x),x)$, with $x\in \beta^k$.

As in the previous section we define a notion of admissibility for the elements of $\mathfrak{S}_{\mathbf{v}}$. Let us emphasize that this notion is different, since we use the notion of spines of a symbol satisfying \eqref{eq:condition}. Fix $S$ a symbol satisfying \eqref{eq:condition}, and denote by $(Y_j)_{j\leq v_1}$ its spines. There are $v_k-v_{k+1}+1$ spines of length $k$ if $k\neq l$ and an infinite number of spines of length $l$. An element $\sigma=(\sigma_j)_{j\leq v_1}\in \mathfrak{S}_{\mathbf{v}}$ will then naturally permute the entries of the spines: if $Y_{i_1}=(\beta_{i_1}^{1},\beta_{i_2}^{2},\ldots,\beta_{i_j}^{j})=(a_1,\ldots,a_j)$ then $Y_{i_1}^{\sigma_{i_1}} =(a_{\sigma_{i_1}{(1)}},\ldots,a_{\sigma_{i_1}{(j)}})$.

We say that $\sigma\in \mathfrak{S}_{\mathbf{v}}$ is {$\heartsuit$-}admissible for the symbol $S$ if
\begin{enumerate}[label=\roman*)]
\item for all $1\leq i \leq l$, the sequence $\gamma^i\in \mathfrak{B}_{v_i}$ obtained by ordering the $i$th entries of all spines of length $\geq i$ is a $\beta$-number, meaning that there are no repeated entries in $\gamma^i$. We then set $S^{\sigma}$ to be the symbol $(\gamma^1,\ldots,\gamma^l)$.
\item $\sigma$ has minimal length among all $w\in \mathfrak{S}_{\mathbf{v}}$ such that $S^{\sigma}=S^{w}$.
\end{enumerate}

Since for $j$ small enough we have $\beta_j^i=j$ for all $1\leq i \leq l$, then the spine {$Y_j$} is $(j,\ldots,j)$. Therefore, any {$\heartsuit$-}admissible $\sigma$ is of finite length and $\sigma_j=\id$ for $j$ small enough. We will denote by ${\mathfrak{S}^{\heartsuit}}(S)$ the set of {$\heartsuit$-}admissible $\sigma\in \mathfrak{S}_{\mathbf{v}}$ for $S$.

\begin{exa}
  We take for $S$ the symbol $S_2$ of the previous example:
  \[
    S = \left(
      \begin{array}{ccccc}
        \ldots & 0 & 2 & 3 & 5 \\
        \ldots & 0 & 2 & 3 & 5 \\
        \ldots & 1 & 3 & 4
      \end{array}
    \right).
  \]
  Then for $j < 0$, the spine {$Y_j$} is $(j,j,j)$ and the spines {$Y_0,Y_1,Y_2$ and $Y_3$} are respectively $(0,0,1)$, $(2,2,4)$, $(3,3,3)$ and $(5,5)$. Note that they do not coincide with the columns.
  
  The element $(\dots,s_2,\id,s_2,\id)$ is then {$\heartsuit$-}admissible and $S^{\sigma}$ is
  \[
    S^{\sigma}= \left(
      \begin{array}{ccccc}
        \ldots & 0 & 2 & 3 & 5 \\
        \ldots & 1 & 3 & 4 & 5 \\
        \ldots & 0 & 2 & 3
      \end{array}
    \right).
  \]
\end{exa}

\begin{Th}
  \label{thm:main_formula}
  Suppose that $S$ is a standard symbol satisfying the condition \eqref{eq:condition}. The canonical basis element associated to the standard $l$-symbol $S$ and multicharge $\mathbf{v}$ is given by
  \[
    G(S) = \sum_{\sigma \in {\mathfrak{S}^{\heartsuit}}(S)}q^{\ell(\sigma)}S^{\sigma}.
  \]
  It is moreover mononial.
\end{Th}

\begin{Rem}
  \label{ord}
  Let $S$ be an ordered symbol satisfying the extra condition of \Cref{ord}. Then the spines are the columns of the symbol, and the two notions of admissiblity coincide. Moreover, the extra condition implies that two consecutive columns do not share a common entry, which implies that for all admissible $\sigma$, we must have $M(S^{\sigma})=0$. Thus \Cref{ordered} and \Cref{thm:main_formula} are consistent with each other.
\end{Rem}

The above theorem will be a consequence of some intermediate results. The strategy is the same as for \Cref{ordered}. We will show that the element $\sum_{\sigma \in {\mathfrak{S}^{\heartsuit}}(S)}q^{\ell(\sigma)}S^{\sigma}$ is monomial and therefore bar invariant. Since it is obviously equal to $S$ modulo $q\mathcal{F}_{\mathbf{v},R}$, this is the canonical basis element associated to $S$.

\begin{lemma}\label{lem:spines}
  Let $S$ be a standard symbol satisfying \eqref{eq:condition}. Assume that $x\in \beta^{i_1}$ and that $x\in \beta^{i_2}$ with $i_1<i_2$. Then for all $i_1<i<i_2$, we have $x\in \beta^i$. 
\end{lemma}

\begin{proof}
  We have $\psi_{i_2,i_1}(x)=x$. For all $i_1<i<i_2$, we then have $\psi_{i_2,i}(x)\leq x$. We thus have $\psi_{i,i_1}(\psi_{i_2,i}(x))=\psi_{i_2,i_1}(x)=x$ thanks to \eqref{eq:condition}. This implies that $x\leq \psi_{i_2,i}(x)$. Therefore $\psi_{i_2,i}(x)=x$ and $x\in \beta^i$.
\end{proof}

Our first goal is to define a monomial element by induction. Let $S$ be a symbol satisfying the condition \eqref{eq:condition} of {nonzero} size. We fix a row $c$ such that all rows below it are the set of $\beta$-numbers of the empty partition: for all $c < i \leq l$, $\beta^i = (j)_{j\leq v_i}$. Such a row exists since $S$ is not of size $0$. Now, let $x\in \mathbb{Z}$ be the smallest such that $x\in \beta^c$ but $x-1\not\in \beta^c$. Such an $x$ exists since the row $\beta^c$ is not equal to $(j)_{j\leq v_c}$. Let $c'\leq c$ be minimal such that $x\in \beta^{c'}$. Then by \Cref{lem:spines}, $x\in \beta^i$ for all $c'\leq i \leq c$.

We have $\beta_{j}^i=j$ for all $i\leq c$ and $j\leq v_{c+1}$. Indeed, for such $i$ and $j$, the standardness of $S$ implies that $\beta_{j}^i \leq \beta_{j}^{c+1} = j$. Since we always have $\beta_{j}^i \geq j$, we obtain the claimed equality. This implies in particular that $x>v_{c+1}+1$: if $j$ is such that $\beta_{j}^c=x$ then $j>v_{c+1}$ and $x>j$. In particular $x,x-1 \not\in \beta^i$ for $i>c$. To summarize, we have
\begin{enumerate}
\item $x\in \beta^c$ and $x-1\not\in \beta^c$,
\item $x\in \beta^i$ for $c'\leq i \leq c$ and $x\not\in \beta^i$ otherwise,
\item for $i>c$, the row $\beta^i$ is $(j)_{j\leq v_i}$ and $v_i<x-1$.
\end{enumerate}

We define a new symbol $S'$ by replacing the values $x$ by $x-1$, but only in the rows $\beta^{c'},\ldots,\beta^c$ which do not already contain $x-1$. The symbol $S'$ is of size smaller than $S$ since, at least the entry $x$ in the row $c$ is replaced by $x-1$.

\begin{lemma}\label{lem:heart}
  The symbol $S'$ constructed above is standard and satisfies \eqref{eq:condition}.
\end{lemma}

\begin{proof}
  We keep the above notations and distinguish three cases:
  \begin{enumerate}
  \item We suppose that $x-1$ appears in the symbol $S$, in a row containing also $x$. Let $d$ be maximal such that $x-1\in \beta^d$, and note that $c' \leq d < c$ since for $i\geq c$, we have $x-1\not\in \beta^i$. Let $d'$ be minimal such that $x-1\in \beta^{d'}$. Then by \Cref{lem:spines}, we have $x-1\in \beta^i$ for all $d'\leq i \leq d$.
    
    The condition \eqref{eq:condition} also implies that $d' \geq c'$. Suppose that $d' < c'$, then $x-1\in \beta^{d'}$ and $x\not\in \beta^{d'}$, which implies that $\psi_{c,d'}(x) = x-1$. But by \eqref{eq:condition}, we have $x-1=\psi_{d,d'}(\psi_{c,d}(x))=\psi_{d,d'}(x)$. Since $x-1$ is in $\beta^d$ and $\beta^{d'}$, we also have $\psi_{d,d'}(x-1)=x-1$, which contradicts the injectivity of $\psi_{d,d'}$.
    
    We also have $c'=1$, that is $x$ is in the first row of the symbol. Suppose that $c'\neq 1$. By choice of $x$, every entry of $\beta^c$ smaller than $x$ is in the rows $\beta^1,\ldots,\beta^{c-1}$, which implies that $y=\psi_{c,1}(x)$ is the biggest entry of $\beta^1$ which is smaller than $x$. The condition \eqref{eq:condition} gives that $y=\psi_{d,1}(x)$, and $\psi_{d,1}(x-1)$ would then be smaller than $y$, contradicting the definition of the map $\psi_{d,1}$.
    
    Therefore, in the symbol $S$, the integers $1 \leq d' \leq d < c$ are such that $x$ is in the row $\beta^i$ if and only if $1\leq i \leq c$ and $x-1$ is in the row $\beta^{i}$ if and only if $d'\leq i \leq d$. The symbol $S'$ is then obtained by replacing $x$ by $x-1$ in the rows $i$ with $1 \leq i < d'$ or $d < i \leq c$. It is almost clear that $S'$ is standard, the only case to consider is the comparison between the rows $d$ and $d+1$ since nothing changed in the row $d$ but the entry $x$ in the row $d+1$ has been replaced by $x-1$, which is a smaller number. The only case where we do not obtain a standard symbol would be if, in $S$, the entries $x$ of rows $d$ and $d+1$ are in the same column:
    \[
      \begin{pmatrix}
        \dots & x-1 & x & \dots\\
        \dots & y & x & \dots
      \end{pmatrix}
      \quad
      \text{is replaced by}
      \quad
      \begin{pmatrix}
        \dots & x-1 & x & \dots\\
        \dots & y & x-1 & \dots
      \end{pmatrix}.
    \]
    By standardness of the symbol $S$, we have $x-1\leq y$ and since $S$ is a symbol, we also have $y < x$, which is impossible. Therefore we deduce that $S'$ is standard.
    
    It remains to show that $S'$ satisfies \eqref{eq:condition}. We will do it by examining how the spines of $S$ behave. Only two spines of $S$ will be modified, the one containing all the entries $x$, and the one containing all the entries $x-1$. Since all the rows below $\beta^c$ are the beta-numbers of the empty partition, the spine of $S$ containing the entries $x$ is of height $c$ and is
    \[
      X=(\underbrace{x,\ldots,x}_{c}).
    \]
    For a similar reason, the spine of $S$ containing all the entries $x-1$ is of height smaller than $c$ and {greater than or equal to} $d$. Let $b$ be maximal such that $x-1$ is in the image of $\psi_{b,d}$ and $z\in \beta^{b}$ is such that $x-1=\psi_{b,d}(z)$, then the spine containing all the entries $x-1$ of $S$ is
    \[
      X_{-}=(\psi_{d',1}(x-1),\ldots,\psi_{d',d'-1}(x-1),\underbrace{x-1,\ldots,x-1}_{d-d'+1},\psi_{b,d+1}(z),\ldots,z).
    \]
    If we denote by $\psi'_{-,-}$ the functions associated to the standard symbol $S'$, for all the values outside these two spines, the functions $\psi_{-,-}$ and $\psi'_{-,-}$ coincide. The spines $X$ and $X_{-}$ are transformed into
    \[
      (\underbrace{x-1,\ldots,x-1}_{d'-1},\underbrace{x,\ldots,x}_{d-d'+1},\underbrace{x-1,\ldots,x-1}_{c-d})
    \]
    and
    \[
      (\psi_{d',1}(x-1),\ldots,\psi_{d',d'-1}(x-1),\underbrace{x-1,\ldots,x-1}_{d-d'+1},\psi_{b,d+1}(z),\ldots,z).
    \]
    These are not spines of $S'$ since the $x$ of the first tuple are exchanged with the $x-1$ of the second tuple to give the following two spines of $S'$
    \[
      X'=(\underbrace{x-1,\ldots,x-1}_{c})
    \]
    and
    \[
      X'_{-}=(\psi_{d',1}(x-1),\ldots,\psi_{d',d'-1}(x-1),\underbrace{x,\ldots,x}_{d-d'+1},\psi_{b,d+1}(z),\ldots,z),
    \]
    which shows that $S'$ satisfies the condition \eqref{eq:condition}.
    
  \item We suppose that $x-1$ appears in the symbol $S$, but no row $S$ contain both $x$ and $x-1$. {The points 2 and 3 before \Cref{lem:heart} imply that} there is no entry $x-1$ in a row below a row containing $x$. Moreover, $x-1$ appears in the row $c'-1$. Indeed, if $y=\psi_{c',c'-1}(x)$, then $y<x$ since $x$ is not in $\beta^{c'-1}$. If $i$ is a row containing $x-1$, then $\psi_{c',i}(x) = x-1$ and $\psi_{c'-1,i}(y) = x-1$ so that $x-1\leq y$. Hence $y=x-1$ and $x-1$ appears in the row $c'-1$. Denote by $k$ the smallest integer such that $x-1$ is in $\beta^k$.
    
    We then replace all $x$ in the symbol $S$ by $x-1$ to obtain $S'$, and it is immediate that $S'$ is standard: if we change an entry $x$ of $S$ into $x-1$, then the entry below in $S$ is either $x$, which is also replaced by $x-1$, {or} greater than $x$ and is still greater than $x-1$.
    
    To show that $S'$ satisfy \eqref{eq:condition}, we do as in the first case. This is easier once we have noticed that all entries $x$ and all entries $x-1$ are in the same spine of $S$. As in the first case, this spine is of height $c$ and is given by
    \[
      X=(\psi_{k,1}(x),\dots,\psi_{k,k-1}(x-1),\underbrace{x-1,\dots,x-1}_{c'-k},\underbrace{x,\dots,x}_{c-c'+1}).
    \]
    The replacement of all $x$ by $x-1$ transforms this spine in
    \[
      X'=(\psi_{k,1}(x),\dots,\psi_{k,k-1}(x-1),\underbrace{x-1,\dots,x-1}_{c'-k},\underbrace{x-1,\dots,x-1}_{c-c'+1}).
    \]
    and it is easily checked that $S'$ satisfy \eqref{eq:condition} since all spines of $S$ are a spine of $S'$ but the spine $X$, which is turned into the spine $X'$.
    
  \item We suppose that $x-1$ does not appear in the symbol $S$. We then replace all $x$ in the symbol $S$ by $x-1$ to obtain $S'$, and it is immediate that $S'$ is standard as in the previous case.

    This is the easiest case for checking the consition \eqref{eq:condition}, since only the spine containing $x$ is modified. Turning all the entries $x$ into $x-1$ will only affect this spine, and we immediately obtaine a spine of $S'$ as in the second case.
  \end{enumerate}
\end{proof}

We can now define a monomial induction on the size. If $S$ is of size $0$, then this is the symbol $\emptyset_{\mathbf{v}}$. Otherwise, consider $x$ and $S'$ as above, and define recusively $\widetilde{\seq}$ by
\[
  \widetilde{\seq}(S) = \widetilde{\seq}(S'),\underbrace{x,\ldots,x}_{m \text{ times}},
\]
where $m$ is the number of $x$ turned into $x-1$ to obtain $S'$. Since $S'$ is standard, satisfies \eqref{eq:condition} and is of smaller size than $S$, the sequence $\widetilde{\seq}(S')$ is well defined.

If $\widetilde{\seq}(S)$ is given by
\[
  \underbrace{x_1,\ldots,x_1}_{\text{$m_1$ times}},\underbrace{x_2,\ldots,x_2}_{\text{$m_2$ times}},\ldots,\underbrace{x_m,\ldots,x_m}_{\text{$m_s$ times}},
\]
we set
\[
  A(S)=F^{(m_s)}_{x_s-1} \ldots F^{(m_1)}_{x_1-1}\cdot\emptyset_{\mathbf{v}}.
\]

\begin{lemma}
  We have
  \[
    A(S)=\sum_{\sigma \in {\mathfrak{S}^{\heartsuit}}(S)}q^{\ell(\sigma)} S^{\sigma}.
  \]
\end{lemma}

\begin{proof}
  We prove this formula by induction on the size of the symbol $S$. The result is clear if $S$ is of size $0$. If $S$ is of {nonzero} size, we choose $x$ and $S'$ as above. Recall that $m$ denotes the number of entries $x$ in $S$ replaced by $x-1$ to obtain $S'$. We also keep all the notations used in the previous lemma and in its proof, since we distinguish the three same cases.
  \begin{enumerate}
  \item The element $x-1$ appears in $S$, in a row containing $x$. We have shown that there are integers $1\leq d' \leq d < c$ such that $x$ is in the row $\beta^i$ if and only if $1 \leq i \leq c$ and $x-1$ is in the row $\beta^i$ if and only if $d' \leq i \leq d$. The $S'$ is obtained by changing the $m=d'+c-d$ entries $x$ into $x-1$ in the rows $\beta^{1},\ldots,\beta^{d'-1},\beta^{d+1},\ldots,\beta^{c}$. Moreover, almost all spines of $S$ is also a spine of $S'$, the only two exceptions being the spines
    \[
      X=(\underbrace{x,\ldots,x}_{c})\quad\text{and}\quad
      X_{-}=(\psi_{d',1}(x-1),\ldots,\psi_{d',d'-1}(x-1),\underbrace{x-1,\ldots,x-1}_{d-d'+1},\psi_{b,d+1}(z),\ldots,z)
    \]
    of $S$, which are turned into the spines
    \[
      X'=(\underbrace{x-1,\ldots,x-1}_{c}) \quad\text{and}\quad X'_{-}=(\psi_{d',1}(x-1),\ldots,\psi_{d',d'-1}(x-1),\underbrace{x,\ldots,x}_{d-d'+1},\psi_{b,d+1}(z),\ldots,z)
    \]
    of $S'$. It is then clear that ${\mathfrak{S}^{\heartsuit}}(S)$ and ${\mathfrak{S}^{\heartsuit}}(S')$ are the same, once we have identified the spines $X$ and $X_{-}$ of $S$ with the spines $X'$ and $X'_{-}$ of $S'$ respectively. Moreover, we have $F_{x-1}^{(m)}\cdot (S')^{\sigma} = S^{\sigma}$ with the above identification of ${\mathfrak{S}^{\heartsuit}}(S)$ and ${\mathfrak{S}^{\heartsuit}}(S')$. Therefore, since by induction hypothesis we have $A(S')=\sum_{\sigma \in {\mathfrak{S}^{\heartsuit}}(S')}q^{\ell(\sigma)} (S')^{\sigma}$, we have $A(S)=\sum_{\sigma \in {\mathfrak{S}^{\heartsuit}}(S)}q^{\ell(\sigma)} S^{\sigma}$.

  \item The element $x-1$ appears in $S$, but not in a row containing $x$. This case is the most difficult one. We have shown that there are integers $1\leq k < c' \leq c$ such that $x$ is in the row $\beta^i$ if and only if $c' \leq i \leq c$ and $x-1$ is in the row $\beta^i$ is and only if $k\leq i < c'$. In this case, $x$ appears $m=c-c'+1$ times in $S$ and $x-1$ appears $c'-k$ times. Therefore, in $S'$, $x-1$ appears $c-k+1$ times and $x$ is not an entry in $S'$. There are then many ways of choosing $m$ entries among the $c-k+1$ entries $x-1$, which is needed to compute $F_{x-1}^{(m)}\cdot (S')^{\sigma}$ for any $\sigma\in {\mathfrak{S}^{\heartsuit}}(S')$.
    
    Fix $\sigma\in {\mathfrak{S}^{\heartsuit}}(S')$, and denote by $1\leq i_1<\dots<i_{c-k+1}\leq c$ the rows containing the entry $x-1$. Then a symbol appears in $F_{x-1}^{(m)}\cdot (S')^{\sigma}$ with a {nonzero} coefficient if and only if it is obtained from $(S')^{\sigma}$ by replacing $m$ entries $x-1$ by $x$. Let $T$ be such a symbol and suppose that the entries $x$ are in the rows $1\leq i_{j_1} < \dots < i_{j_m} \leq c$. Then, using \Cref{divided}, the coefficient of $T$ in $F_{x-1}^{(m)}\cdot (S')^{\sigma}$ is $q^{\sum_{r=1}^{m}(c-k+1-j_r-m+r)}$ since, below the entry $x$ in the row $i_{j_r}$, there are $m-r$ entries $x$ and $c-k+1-j_r-m+r$ entries $x-1$.

    We now compare ${\mathfrak{S}^{\heartsuit}}(S)$ and ${\mathfrak{S}^{\heartsuit}}(S')$. Only the spine
    \[
      X=(\psi_{k,1}(x),\dots,\psi_{k,k-1}(x-1),\underbrace{x-1,\dots,x-1}_{c'-k},\underbrace{x,\dots,x}_{c-c'+1})
    \]
    of $S$ is not a spine of $S'$, and is replaced by the spine
    \[
      X'=(\psi_{k,1}(x),\dots,\psi_{k,k-1}(x-1),\underbrace{x-1,\dots,x-1}_{c'-k},\underbrace{x-1,\dots,x-1}_{c-c'+1})
    \]
    of $S'$. Hence every {$\heartsuit$-}admissible permutation of the spine $X$ will be obtained from first applying a {$\heartsuit$-}admissible permutation of the $c-k+1$ last entries of the spines (we permute the various $x$ and the $x-1$) and then a {$\heartsuit$-}admissible permutation of the spine $X'$ of $S'$. Therefore, every $\sigma\in {\mathfrak{S}^{\heartsuit}}(S')$ gives rise to $\binom{c-k+1}{m}$ elements of ${\mathfrak{S}^{\heartsuit}}(S)$ by precomposing by permutations of the spine $X$ permuting the last $c-k+1$ entries. The symbol $T$ above, is then obtained from $S$ by precomposing $\sigma$ with a permutation of $X$ moving the $m$ entries $x$ at the positions $j_1+k-1< j_2+k-1 < \dots < j_{m}+k-1$, and the minimum length of such a permutation is $\sum_{r=1}^m(c-m+r-(j_r+k-1))$. This is exactly the power of $q$ in the coefficient of $T$ in $F_{x-1}^{(m)}\cdot (S')^{\sigma}$.
    
    Therefore, using the induction hypothesis, we obtain that $F^{(m)}_{x-1}\cdot A(S') = \sum_{\sigma\in{\mathfrak{S}^{\heartsuit}}(S)}q^{\ell(\sigma)}S^{\sigma}$ as expected.
    
  \item The element $x-1$ does not appear in $S$. All spines of $S$ are spines of $S'$, with the exception on the one containing $x$ which is turned into the spine of $S'$ containing $x-1$. The identification between ${\mathfrak{S}^{\heartsuit}}(S)$ and ${\mathfrak{S}^{\heartsuit}}(S')$ is clear, as well as the equality $F_{x-1}^{(m)}\cdot (S')^{\sigma} = S^{\sigma}$. Therefore, since by induction hypothesis we have $A(S')=\sum_{\sigma \in {\mathfrak{S}^{\heartsuit}}(S')}q^{\ell(\sigma)} (S')^{\sigma}$, we have $A(S)=\sum_{\sigma \in {\mathfrak{S}^{\heartsuit}}(S)}q^{\ell(\sigma)} S^{\sigma}$.
  \end{enumerate}
\end{proof}

\begin{proof}[Proof of \Cref{thm:main_formula}]
  The element $A(S)$ is monomial, hence bar invariant. Moreover, it is clear that $A(S) = S \mod q\mathcal{F}_{\mathbf{v},R}$ since $\ell(\sigma)\geq 0$ with equality if and only if each component of $\sigma$ is the identity. Therefore, $A(S)$ is the canonical basis element $G(S)$.
\end{proof}

\section{Additional remarks}

Most of our results are available here when ${\bf v} \in \mathcal{A}^l$. One may ask what can be expected in the general case. This is particularly important when one wishes to apply these results in the context of graded representation theory, as described by Brundan and Kleshchev \cite{BK}.

The symmetric group $\mathfrak{S}_l$ acts naturally on $\mathbb{Z}^l$. If $\mathbf{v} = (v_{1},\ldots,v_{l})\in\mathbb{Z}^{l}$, there exist $\mathbf{s} = (s_{1},\ldots,s_{l})\in\mathcal{A}^{l}$ and $w\in \mathfrak{S}_l$ such that $\mathbf{v}=w\cdot\mathbf{s}$. Since $\mathbf{v}$ and $\mathbf{s}$ yield the same dominant weight, we have an isomorphism $\phi_{\mathbf{s},\mathbf{v}}$ of $U_q(\mathfrak{gl}_{\infty})$ from $V_{\bf s}$ to $V_{\bf v}$. Similarly, the canonical basis of the module $V_{\bf v}$ is naturally indexed by a certain subset of symbols $\mathcal{S}_{\bf v}$ and, by \cite{JL}, we have a bijection $\varphi_{\mathbf{s},\mathbf{v}} : \mathcal{S}_{\bf s}^{st} \to \mathcal{S}_{\bf v}$ which can be combinatorially described. 

Finally, by the uniqueness of the crystal basis proved by Kashiwara, we also have:
\[
\phi_{\mathbf{s},\mathbf{v}}(G (S) ) = G(\varphi_{\mathbf{s},\mathbf{v}}(S)).
\]
Thus, our results can be applied to obtain similar formulas in the general case (in particular, the canonical basis retains the ``monomial'' property under this isomorphism).

We also remark that recently, Lyle \cite{L1,L2} has obtained a general closed formula for the canonical basis elements for arbitrary choices of ${\bf s}$. However, Lyle's results seem to differ from ours since she works with a $U_q(\widehat{\mathfrak{sl}}_{e})$ version of the Fock space. For instance, in \cite{L2}, these results concern the case where the decomposition numbers are indexed by two symbols having the same ``multicore''. In contrast, in our case with $e=\infty$, every symbol is a multicore and Lyle's result only implies that the coefficient of $S$ in $G(S)$ is $1$, which is part of the definition of $G(S)$.

\noindent {\bf Addresses} : \\N.J :Laboratoire de Math\'ematiques de Reims (UMR CNRS 9008), Universit\'{e} de Reims Champagne-Ardennes, UFR Sciences exactes et
naturelles.  Moulin de la Housse BP
1039. 51100 Reims. France. nicolas.jacon@univ-reims.fr \\
\\
A.L: Laboratoire de Math\'ematiques Blaise Pascal (UMR CNRS 6620),
Universit\'e Clermont Auvergne,
Campus Universitaire des C\'ezeaux,
3 place Vasarely,
TSA 60026
CS 60026
63178 Aubi\`ere Cedex - France. abel.lacabanne@uca.fr

\end{document}